\documentclass{article}
\usepackage{graphicx}
\usepackage {color}
\font \m=msbm10
\newcommand{\R}{{\hbox {\m R}}}
\newcommand{\C}{{\hbox {\m C}}}

\newcommand{\U}{{\hbox {\m U}}}
\newcommand {\HH}{{\hbox {\m H}}}
\newcommand {\eps} {{\varepsilon}}
\newtheorem {lemma}{Lemma}
\newcommand {\lu}{{\lambda}}

\newtheorem {theorem}{Theorem}
\begin{document}
\author{Wendelin Werner}
\title{Some recent aspects of random conformally invariant systems}
%
\date {Universit\'e Paris-Sud and Ecole Normale Sup\'erieure}
\maketitle 

\begin {abstract}
These are the lecture notes from a course given in July 2005 at the summer school in 
Les Houches. We describe some recent results concerning two-dimensional 
conformally invariant systems. In particular, 
we discuss conformally invariant measures
on loops and conformal loop-ensembles (CLE). 
\end {abstract}   
%

\eject

\tableofcontents

\eject

\section*{Overview}

The goal of these lectures is to give a self-contained account of recent developments concerning the theory of 
continuous two-dimensional and conformally invariant structures. These object are conjectured -- and in some cases it is proved -- to arise as scaling limits of various critical models from statistical physics in the plane. These include percolation, the Ising model, the Potts models, the $O(N)$ models, self-avoiding walks, loop-erased and Hamiltonian random walks, the uniform spanning trees. We shall briefly define these models in our first lecture. 

I will try to put the emphasis on one ongoing line of research
rather than repeating some earlier courses and lectures on the subject. In particular, I will neither
spend much time describing the definition of Schramm-Loewner Evolutions (SLE)
nor discuss details of the aspects of the theory that have been developed prior to 2002 as this material is already available in the form of lecture notes (including by myself) \cite {Lln, Wln, Wln2}, survey papers for theoretical physicists \cite {KNln, Cln},
and now even in an excellent book \cite {Lbook}.  The list of references of the present notes are therefore by no means complete (I just listed the papers that are directly relevant to the content of the notes, as well as some 
historical  or survey-type papers) and it should not be considered as a ``selection'' of recent papers of the subject. One could mention many nice recent papers by M. Bauer, R. Bauer, Bernard, Dub\'edat, Friedrich, Makarov, Schramm, Sheffield, Wilson, Zhan etc. that are not listed here. 

Sections one and three of the present lectures/notes are going to be more classical (discrete models, SLE definition), while the second and forth one will present recent or ongoing research.

In the first one, I will mostly present some discrete lattice-based 
models originating in the physics literature,
 and state precisely the conformal invariance conjectures. They roughly say that when
 the mesh-size of the lattice vanishes, the models have  (scaling) limits, and that
 these limits are conformally invariant.
It is important to recall that in most cases, but with certain important exceptions, these conjectures are not proved (even on heuristic level), and that they rely mostly on 
numerical evidence, and (now) on the existence of nice continuous models that they would converge to. One of the ``fronts'' on which mathematicians are working is to try to prove these conjectures and one can hope for progress in the future. 
The problem is to work on specific discrete 
models, to understand some of their deep properties that enable to control their large-scale behavior.  
This important 
question will however not be the main focus of these lectures.

Instead, motivated by the quest for the understanding of the large-scale properties and scaling limits of these models, we shall present some random continuous objects (random loops, random curves, families of random loops) and state
various results that show that IF the conformal invariance conjecture hold, then these continuous models are the only possible scaling limits of the discrete models.

In particular, in the second lecture, I will construct an important (infinite) measure supported on 
the set of self-avoiding loops in the plane, 
and show that it is the unique conformally invariant measure that satisfies a so-called restriction property (that is shared in the discrete case by Brownian loops boundaries, percolation cluster boundaries, self-avoiding loops).  This will show that the scaling limits of these discrete models ought to be the same.  
The restriction property
 can be interpreted as an invariance under the action of a certain semi-group of
 conformal maps.
The results presented in this lecture are derived in \cite {Wsal} and 
inspired by ideas developed with Greg Lawler and Oded Schramm in \cite {LW2, LSWrest, LSWsaw}.

In the third lecture, I will focus on the law of one curve (for instance an interface), and the ideas behind the definition of SLE via iterations of random conformal maps, and survey some of the consequences. This will be reminiscent of the above-mentioned more detailed papers/lecture notes/books on SLE. This can serve as an introduction to SLE for those who are not acquainted to it. The readers who are may want to skip this third section.

In the fourth lecture, I shall focus on random continuous families of non-overlapping loops that possess a certain (natural in view of the discrete models) conformal invariance property introduced in \cite{ShW}.
 These are (now) called 
the ``Conformal loop-ensembles'' (CLE in short). 
I shall define this property and discuss some of its first consequences.
This will lead rather naturally to the definition of the Brownian loop-soup and to the idea that the outer boundaries of ``subcritical'' loop-soups provide examples of CLEs (and in fact these are all the possible CLEs).
Note that this part is ``almost SLE independent''.
We will conclude these lectures by presenting briefly the Gaussian Free Field (GFF in short) and some of its properties. Scott Sheffield, together with Oded Schramm \cite {ScSh, Sh}, is currently studying the ``geometry'' of this generalized random function in terms of SLE curves, and it turns out that  one can find (in a deterministic way) CLEs inside the Gaussian Free Field, and that conversely, the field can be recovered by the CLE. Hence, there is a link between GFF, CLEs, loop-soups and SLEs. 

\eject

\section {Some discrete models}

In this section, we shall describe quickly
 various discrete systems that are conjectured to be conformally invariant in the scaling limit. All these discrete systems have been introduced in the physics 
literature.
They live on (finite or infinite portions) of a planar lattice. 
We shall not discuss here the rather delicate issue of for what lattices the conformal invariance conjecture should hold. So, to fix ideas, we will assume that the lattice that we use is the square lattice, or the triangular lattice, or its dual -- the honeycomb lattice.
We will in fact define only two types of measures:
\begin {itemize}
\item
Infinite measures on the set of connected subsets of the lattice (for instance polygons). These measures will be 
translation-invariant (and therefore clearly infinite).
\item
Probability measures on ``configurations'' in a given finite portion of the lattice
(i.e. in a fiven domain). 
\end {itemize}
We will first describe each of these models separately, and then, in the final subsection, we shall state a precise version of the (sometimes proved) conformal invariance 
conjecture.

\subsection {Self-avoiding walks and polygons}

We say that a nearest-neighbor path of length $N$ started from the origin, and with no double point is a self-avoiding path of length $N$ (on the chosen lattice). The number of such walks will be denoted by $A_N$. 
It is easy to see that $A_{N+M} \le A_N A_M$ and that $A_N \ge c^N$ for some $c>1$ (by counting the number of ``north-east'' paths that are anyway self-avoiding), and this implies immediately that  there exists $\lu >1 $ such that $\lim_{N \to \infty} (A_N)^{1/N} = \lu$. In other words, the number 
of self-avoiding walks of length $N$ on the considered lattice behaves roughly like $\lu^N$. It is easy to see that 
$\lu$ is a lattice-dependent constant.

\begin {figure}[b]
\begin {center}
\begin{picture}(0,0)%
\includegraphics{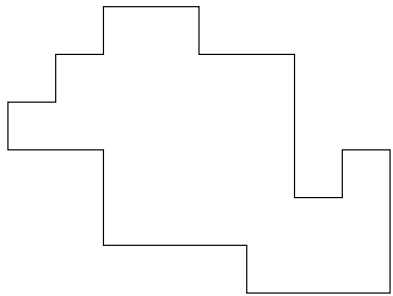}%
\end{picture}%
\setlength{\unitlength}{3947sp}%
\begingroup\makeatletter\ifx\SetFigFont\undefined%
\gdef\SetFigFont#1#2#3#4#5{%
  \reset@font\fontsize{#1}{#2pt}%
  \fontfamily{#3}\fontseries{#4}\fontshape{#5}%
  \selectfont}%
\fi\endgroup%
\begin{picture}(2424,1824)(1489,-2773)
\end{picture}%

\caption {A self-avoiding polygon on the square lattice}
\end {center}
\end {figure}

The length of a walk is not a quantity that is well-suited to conformal invariance. One can therefore 
look for self-avoiding walks with prescribed end-points but varying length, as described for instance 
in \cite {LSWsaw}. 
Another interesting 
possibility is to consider self-avoiding polygons instead: We say that a nearest-neighbor path $(w_0, \ldots w_N)$ of length $N$ is a self-avoiding polygon of length $N$ if
$w_0=w_N$ and if it is otherwise self-avoiding.
Again, it is not difficult to see that for the same constant $\mu$ as before, the number $A_N'$ of self-avoiding polygons of length $N$ with $w_0=0$ behaves like $\lu^{N + o(1)}$ as $N \to \infty$.

We are from now on going to identify two polygons $w$ and $w'$ of length $N$ if for some $n_0$, 
$w_{n_0+n} = w_{n}'$ for all $n$ (and $w$ is extended into a $N$-periodic function). In other words, there is no 
special marked point (i.e. no root) on a self-avoiding polygon.
Hence, $A_N'$ is the number of self-avoiding polygons of length $N$ that pass through the origin.

We define the measure $\mu_{sap}$ on the set of self-avoiding polygons on the lattice in such way 
that for any given self-avoiding polygon $w$ of length $N=N(w)$, 
$$ 
\mu_{sap} (\{ w \}) = \lu^{-N}.
$$ 
This is clearly an translation-invariant infinite measure. 
Note that in order to define this measure, one does not need to specify any marked points or any prescribed domain.

\subsection {Random walks loops}

Suppose that each point on the lattice has $d$ neighbors (i.e. $d=3, 4, 6$ respectively on the hexagonal, square and 
triangular lattices). A loop is a nearest-neighbor path 
$w=(w_0, w_1, \ldots, w_n= w_0)$ of a certain length $n$, 
such that $w_0=w_n$, and as before, we identify two loops modulo shift of the time-parametrization. 
The difference with the self-avoiding loops is that $w$ is allowed to have multiple points.

Now, define the measure $\mu_{rwl}$ on the set of loops
that puts a weight $d^{-n}$ to each loop of length $n$ (for all $n>0$). 
Again, this is a translation-invariant measure.

\subsection {Site-Percolation}

For each site of the considered
infinite lattice, toss a coin independently: With probability $p$ the site is declared open, and with probability $1-p$, it is declared vacant. Then, one is interested in the connectivity properties of the set of open sites. 
One can consider the ``connected components'' of open sites. A connected component $C$ is a connected family of open sites, such that all the sites neighboring $C$ are closed.

For planar graphs, and at least for nice lattices, such as the ones that we shall concentrate on, 
one can prove (e.g. \cite {G}) that there exists a critical value $p_c$ such that:
\begin {itemize}
\item If $p \le p_c$, then with probability one, there is no infinite connected component of open sites.
\item If $p > p_c$, then with probability one, there is a unique infinite connected component of open sites.
\end {itemize}
We will from now focus solely on the critical case where $p=p_c$. Then, one can prove (via the so-called Russo-Seymour-Welsh arguments that are by the way not established for all planar lattices -- one needs symmetry conditions, but these conditions are 
however fulfilled in the cases that we are considering, see \cite {Kbook}) that the probability of existence of open connected components
at any scale is bounded away from zero and from one. More precisely, if $A$ and $B$ are two non-empty open sets (at positive distance from each other) in the plane, then, there exists $\eps >0$ such that for any sufficiently large $N$, the probability $P(NA,NB)$ that there exists a connected component of open sites that intersects both $NA$ and $NB$ satisfies $\eps < P(NA,NB) \le 1 - \eps$.

The percolation probability measure defines a random collection of disjoint clusters of open sites $(C_j, j \in J)$.
All clusters are finite, but there exist clusters at any scale because of the Russo-Seymour-Welsh type arguments. 

Define the infinite measure  $\mu_{pcl}$ on the set of 
possible clusters 
such that 
$$ \mu_{pcl} ( \{ C \} ) = P ( C \in \{ C_j, j \in J\}).
$$
In other words, for any family ${\cal C}$ of clusters,  then $\mu_{pcl} ({\cal C}) $ is the expected number of percolation clusters that occur and belong to ${\cal C}$ for one realization.
Another possible direct description of $\mu_{pcl}$ is 
$$
\mu_{pcl}( \{C \})= p_c^{\# C} (1-p_c)^{\# \partial C},
$$
where $\partial C$ is the set of sites that are at unit distance of $C$.

One can equivalently describe a cluster via the collection of its ``boundaries''. For instance, its 
outer boundary is a path that goes through all points of the cluster that are accessible from infinity.

The value of the critical probability is lattice-dependent. 
The triangular lattice is special: A simple symmetry argument shows intuitively (and this can be made rigorous, e.g. \cite {Kbook}) that for this lattice $p_c= 1/2$. 
In this case, the outer boundary of a cluster can be viewed as a self-avoiding loop drawn on a
hexagonal grid.
The measure $\mu_{pcl}$ defines a measure on outer boudaries of clusters that we call $\mu_{perco}$.
The $\mu_{perco}$ mass of a given self-avoiding loop
 is the probability that it is the outer boundary of a  
cluster. Hence, for a self-avoiding loop $w$ on the hexagonal lattice,
 if $k(w)$ is the number of its neighboring hexagons,
$$\mu_{perco} ( l ) = 2^{-k(w)}.$$
 
It is also possible to define another version of percolation, where this time, the edges of the graph are independently declared to be open or closed with probability $p$. One can then also define a critical probability $p_c$, clusters of open sites and their boundaries. In this case, it turns out that the square lattice is special in that $p_c = 1/2$.
 
\subsection {The Ising model}

We are now going to consider a system similar to the site percolation described above.
In the statistical physics language, this is the Ising model without external magnetic field -- and everybody in this audience probably knows it.
To define it, it is convenient to focus on a finite portion of the infinite lattice. One can think of this graph as the intersection of a (large but bounded) open subset $D$ of the plane with the lattice. 

Each site $x$ can take one of two colors (or spins): $\sigma (x)$ can be equal to $1$ or to $2$ (a configuration $\sigma$ is therefore a function that assigns to each site $x$ a color $\sigma(x)$). 
Two neighbors prefer however to have the same color: For each configuration $\sigma$, we define the Hamiltonian 
$$ H(\sigma)= \sum_{x \equiv y}  1_{\sigma (x) \not= \sigma (y)} $$
where the sum is over all pairs of neighboring sites $x$ and $y$ in the (finite) considered graph. 

For each value of the positive parameter $\beta$, we are now going to define a probability measure $P_\beta$  on 
the set of configurations: 
$$
P_\beta (\sigma) = \frac {1}{Z_\beta} \exp \{- \beta H (\sigma) \},
$$
where the partition function $Z_\beta$ is defined to be 
$ Z_\beta = \sum_\sigma \exp \{- \beta H(\sigma) \}$ so that $P_\beta$ is indeed a probability measure.

When $\beta$ is small, the interaction between spins at neighboring sites is small, but if $\beta$ is large, it is strong. In fact, when $\beta=0$ then one just has percolation with $p=1/2$, and in the  $\beta \to \infty$
limit, all sites have almost surely the same color.

It turns out that for very large graphs (and this can be rigorously stated in different ways -- for instance in the limit when the graph tends to fill the infinite lattice), one observes different behaviors depending on the 
value of the parameter $\beta$. 

One can note at this point (this is for instance clear from the relation with percolation at $p=1/2$ for very small $\beta$) that for general lattices, the existence or not of very large (infinite) clusters of sites of the same color (we shall say monochromatic clusters) is not necessarily the suitable or natural observable to study.
 We shall however come back to this in the context of the $O(N)$ models, where one sticks to the triangular lattice.

Physically, for the Ising model, the most relevant quantity is not related to the connectivity property of the 
configuration, but rather to the ``correlation functions''. The simplest such function is 
the two-point function, which is simply the correlation between the value of the spin at two different sites $x$ and $y$:
$$ 
F_\beta (x,y) = P_\beta ( \sigma(x) = \sigma (y) ) - 1/2
.$$
Since $\sigma(x)$ and $\sigma (y)$ are easily shown to be ``positively correlated'', one has that $F > 0$.

When $\beta$ is small, then the correlation function tends exponentially fast to zero (when the sites are fixed and the mesh-size of the lattice goes to zero -- or when the lattice is fixed and the distance between the sites increases).
When $\beta$ is too large, then in fact $F$ is bounded away from zero. This corresponds intuitively to the idea that one of the two colors ``wins'' and occupies  more sites than the other.
The critical value of $\beta$ is the one that separates these two phases. In fact, it can be shown on the square lattice (this is Onsager's calculation) that at this critical value, 
$F$ decays like a certain power law in terms of the mesh-size of the lattice. 
 
\subsection {The Potts models}

The Potts model is a natural generalization of the Ising model: This time, a site can choose between a given number 
$q \ge 2$ of colors. The Hamiltonian and the probability measures $P_\beta$ are then defined exactly in the same way as for the Ising model (and when $q=2$, the Potts model is exactly the Ising model). 

One is then interested in the correlation functions 
$$
 F_\beta (x,y) = P_\beta ( \sigma (x) = \sigma (y) ) - 1/q .
$$ 
Again, one observes a supercritical phase (where on color ``wins'') and a subcritical phase (with exponentially decaying correlations). However, when $q$ is large, it is believed (in fact it is proved for very large $q$) that at the critical value of $\beta$, one still has one ``winning color''. Hence, the critical model will exhibit interesting phase coexistence properties only for small values of $q$. The conjecture is that in two dimensions, this will hold up to $q=4$.

For both the Ising model and the Potts model, we have just described the model with so-called ``free'' boundary conditions. It is possible to condition the probability measure by constraining the sites on the boundary of the domain to take a given color or given colors (one can divide the boundary into two parts and assign a color to each part). At the critical value, and on large scale, these boundary conditions are important and their influence does not really ``disappear'' (of course this notion depends on what one is looking at...) in the scaling limit.

\subsection {FK representations of Potts models}

We just mentioned that the notion of ``Ising clusters'' or ``Potts model clusters'' was not so natural when the lattice is not the triangular lattice. One can however relate properties of the Potts models to connectivity properties of another model, which is a dependent percolation model, introduced by Fortuin and Kasteleyn in the 60's, 
 also sometimes known as the random-cluster model. This is a probability measure on configuration on ``edges'' of a graph: Each edge of a finite graph can be declared open or closed. So, each configuration $\omega$ has a number 
$o(\omega)$ of open edges, a number $c(\omega)$ of closed edges, and it defines $k(\omega)$ different 
open clusters (i.e. sets of sites connected by open edges). 

The FK percolation measure with parameters $p \in (0,1)$ and $q > 0$ is defined as 
$$
P( \omega) = \frac 1 {Z_{p,q}} p^{o(\omega)} (1-p)^{c(\omega)} q^{k(\omega)},
$$
where as before, $Z$ is the normalizing constant.

When $q=1$, this is just the percolation on edges mentioned before. In fact, when $q$ is an integer, this FK measure is closely related to the (free) Potts model with the same $q$. It is indeed easy to check that if one first chooses 
$\omega$ according to the FK measure, then colors each of the $k(\omega)$ clusters independently by one of the $q$
colors, then the obtained coloring has exactly the law of a Potts model, for a given $\beta (p)$. In particular, this \begin {eqnarray*}
F(x,y) = P ( \sigma (x) = \sigma (y)) - 1/q & =& P ( x \leftrightarrow y ) + P (x \not\leftrightarrow y)/q
- 1/q \\
&=&  (1 - 1/q) P( x \leftrightarrow y) 
\end {eqnarray*}
where the symbol $\leftrightarrow$  refers to the connectivity properties of $\omega$ and $F$ is the two-point function of the 
Potts model. Hence, when $q$ is fixed, one can relate the phase-transition of the Potts model to a phase-transition in terms of connectivity properties of the random cluster measure (letting $p$ vary) for very large graphs. Note however that while the Potts models made sense only for integer values of $q$, this FK percolation allows any positive real value of $q$.

The conjecture is that in two dimensions, as for the Potts models,  at the critical value $p_c$, the measure exhibits nice large-scale properties up to $q=4$ only.

When $q \to 0+$ and $p=p_c (q)$, the model converges to the ``uniform spanning tree''. This is the uniform distribution on the set of configurations with one single cluster (all sites are connected) but no open cycle (i.e. the graph is a tree). This is a combinatorially rich model, that can be constructed using ``loop-erased random walks'' or domino tilings (there is a bijection between the set of domino tilings and the set of uniform spanning trees) for which other tools (explicit determinant computations for correlations) are also available.

An interesting property of the FK model in the plane is its ``duality''. If one looks at the dual graph
$\omega^*$ of a sample $\omega$ of an FK model, then $\omega^*$ also follows an FK probability measure, on the 
dual graph. In particular, this leads to the notion of ``wired'' boundary conditions, where all vertices on the 
exterior boundary of a finite planar graph are identified, see e.g. \cite {G2}:
The dual of the ``free'' FK model is the ``wired'' one.

\subsection {The $O(N)$ models}

This is a model that makes use of specific properties of hexagonal/triangular 
lattices. We shall be focus here on the 
honeycomb lattice, where each site has exactly three incoming edges, and each face has six sides, which is the dual of the regular triangular lattice.
 A configuration $c$ on (a portion) of the lattice will be a family of disjoint finite self-avoiding loops on this 
 lattice. For each configuration, there are $l(c)$ loops, and the total cumulated length of the loops is 
 $L(c)$ (counted in number of steps). 
 We assign to each configuration $c$, a probability proportional to 
 $N^{l(c)} \theta^{L(c)}$, where $N >0 $ and $\theta < 1$ are two parameters of the model.
The partition function is then 
$$Z = \sum_c N^{l(c)} \theta^{L(c)}$$
where the sum is over all admissible configurations (in general all possible configurations in a given finite 
portion of the lattice).
 
\begin {figure}
\begin {center}
\includegraphics [width=1.5in]{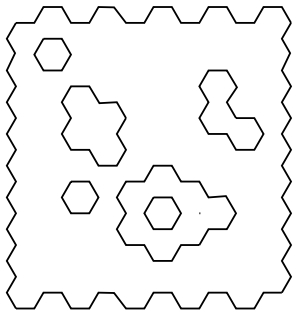}
\caption {An $O(N)$ configuration}
\end {center}
\end {figure}

It is worthwhile noticing that when $N=1$, one can view the loops as the phase separation lines of the 
Ising model on a triangular lattice, with $\theta$ related to the parameter $\beta$ of the Ising model 
by $\exp (- \beta) = \theta$. 
Also, when $N \to 0$, one can relate a conditioned version of the $O(N)$ measure, to some features of the measure on self-avoiding loops (just condition on the existence of at least one loop and let $N \to 0$).
For general integer $N$, the $O(N)$ model is reminiscent of the $N+1$-Potts model on the triangular lattice, but there is the important additional constraint that no 
small triangle has corners of three different colors.

It is conjectured that when $N$ is not too large and fixed ($N \le 2$ in fact), the $O(N)$ model undergoes a phase transition when $\theta$ varies. And, at a critical value $\theta_c(N)$, one observes an interesting
large-scale phenomena (actually the phase where $\theta > \theta_c $ is interesting too), see e.g. \cite {N,KNln}. 
Nienhuis conjectured that  $\theta_c (N)$ is equal to $1/\sqrt { 2 + \sqrt {2-N}}$.
In patricular, the constant $\lambda$ describing the number of self-avoiding walks on that lattice should be $1/ \theta(0+)= \sqrt {2 + \sqrt {2}}$.

Hence, for all $N \le 2$, one has discrete probability distributions on families of disjoint self-avoiding loops in a domain, that are conjectured to behave ``nicely'' in the scaling limit.

\subsection {Conformal invariance}

Let us first concentrate on the ``probability measures'' on configurations. 
For all previous critical models, let us consider a fixed open subset $D$ of the complex plane. For each small $\delta$, one considers a suitable discrete approximation of $D$ by a subset of the graph $\delta {\cal G}$. We 
call $D_\delta$ this graph.

We have for each $\delta$
a probability measure on configurations of disjoint loops (for the $O(N)$ models) or clusters (for the 
FK percolation) in $D_\delta$.

\medbreak
{The {\em conformal invariance conjecture} for these probability measures
 is that (for each $D$):
 \begin {itemize}
 \item {(i)} These probability measures converge when $\delta \to 0$,
toward a probability measure $P_D$ on ``continuous configurations'' in the domain $D$.
One plausible notion of convergence is for instance, for each $\eps >0$, the convergence of any continuous (in the Hausdorff metric) functional that depends only of those clusters/loops of diameter greater than $\eps$.
 \item {(ii)} If $\Phi$ is a conformal map defined on $D$, then the law of the image of $P_D$ under the conformal map $\Phi$ is identical to $P_{\Phi (D)}$.
\end {itemize}}

\medbreak
 
At this point, it is time to recall a few definitions and facts about conformal maps, and in particular Riemann's mapping theorem.
\begin {itemize}
\item If $D$ and $D'$ are two simply connected domains in the plane, we say that $\Phi : D \to D'$ is a 
conformal map from $D$ onto $D'$ if it is an angle-preserving bijection from $D$ onto $D'$ (angle preserving for instance means that the images of two smooth curves that intersect with angle $\alpha$ are two smooth curves that intersect  with angle $\alpha$ -- it is easy to see that this is equivalent to say that at each point, $\Phi$ is 
differentiable in the complex sense).
\item There exists a three-dimensional family of conformal maps from the unit disc onto itself: These are the Moebius transformations, that one can explicitly describe and for which one can directly check that they are conformal maps.
It is also easy to check that they are no other such maps.
\item For any two simply connected open proper 
subsets $D$ and $D'$ of the complex plane, one can find angle-preserving bijections (i.e. conformal maps) from
$D$ onto $D'$. In fact, there exists a three-dimensional family of such maps (this follows from the previous item).
This is Riemann's mapping theorem. Its proof can be found in any complex analysis textbook.
\end {itemize}
\begin {figure}[t]
\begin {center}
\begin{picture}(0,0)%
\includegraphics{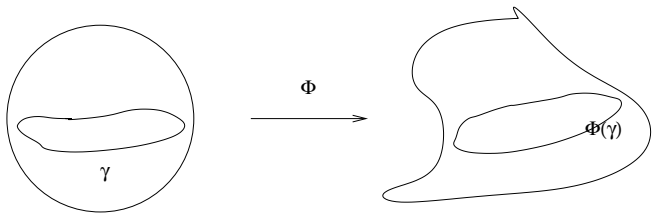}%
\end{picture}%
\setlength{\unitlength}{1381sp}%
\begingroup\makeatletter\ifx\SetFigFont\undefined%
\gdef\SetFigFont#1#2#3#4#5{%
  \reset@font\fontsize{#1}{#2pt}%
  \fontfamily{#3}\fontseries{#4}\fontshape{#5}%
  \selectfont}%
\fi\endgroup%
\begin{picture}(11589,3715)(1009,-3853)
\end{picture}%

\caption {The image of a loop under a conformal map (sketch).}
\end {center}
\end {figure}

 Hence, the conformal invariance conjecture implies in particular that the probability measure $P_D$ is 
 invariant under the group of conformal maps from $D$ onto itself.
 
 \medbreak
 
In the case of the infinite translation-invariant measures on loops, clusters or self-avoiding loops, the conformal invariance conjecture goes similarly: Consider the measures on the lattice with mesh-size $\delta$, and define the measures $\mu_{D,\delta}$ as the restrictions of these measures to those loops/clusters that stay in $D_\delta$. Then, these measures are conjectured to converge weakly to a measure $\mu_D$ on subsets of $D$, and this limiting measure should be conformally invariant.
 
 \medbreak
 
Among these conjectures, only the following have at present been proved:

\begin {itemize}
\item
The measure on random walk loops converges to a continuous measure on Brownian loops, and it follows from the 
conformal invariance of planar Brownian motion -- known since Paul L\'evy -- that the measure on Brownian loops is conformally invariant \cite {LWls}. See also \cite {LT} for convergence issues.

\item
Critical percolation on the triangular lattice (recall that this is the special case where $p_c=1/2$) has been 
proved by Smirnov \cite {Sm} to be conformally invariant. 
This implies the conformal invariance conjecture for the measure $\mu_{perc}$,
see \cite {CN}. 

\item
The $q=0+$ limit of the FK percolation model (the uniform spanning tree) has been proved \cite {LSWlesl} to 
converge to a conformally invariant object (here one needs to keep track of the paths in the tree, the Hausdorff topology is not well-suited).

\end {itemize}

\medbreak

In particular, it is still a challenge to prove conformal invariance of self-avoiding walks, percolation on other lattices, FK percolation and Ising/Potts models.

\medbreak

Let us at this point make a general comment on the theory of conformal maps: A conformal map can be viewed 
as acting globally on a domain, and this is the point of vue that we are going to emphasize in these lectures.
On the other hand, a conformal map $\Phi$ defined on $D$ is holomorphic at every point $z_0$ of $D$,  and it admits a series expansion $\Phi (z_0+h) = \sum_{n\ge 0} a_n h^n$ in the neighborhood of $z_0$. This local description of $\Phi$ can also be very useful for instance, if one composes conformal maps. Very very loosely speaking, ``conformal field theory'' that had been successfully developed by theoretical physicist to understand the scaling limit of 
the above-mentioned critical models, uses more this local aspect of conformal maps.

\medbreak

Also, it is important to mention that a closely related issue is the ``universality conjecture'': If one considers the same model (say, the critical Potts model for $q=3$) on two different lattices (say the square and the triangular), then the limiting continuous measure should be the same in the two cases, despite the fact that the critical value of the parameter $p$ (resp. $\theta$ or $\beta$) is lattice-dependent. As we shall see in the current lectures, this can in fact be viewed as a (non-trivial) consequence of the conformal invariance conjecture, since we will show that only a limited number of conformal invariant possible scaling limits do exist.

\eject
 
\section {A ``conformal Haar measure'' on self-avoiding loops}

\subsection {Preliminaries}

A standard way to state Riemann's mapping theorem goes as follows: If $D$ and $D'$ are two imply connected subsets $\C$, that are not equal to $\C$, and that contain the origin, then there exists unique conformal map 
from $D$ onto $D'$ such that $\Phi (0)= 0$, and $\Phi'(0)$ is a positive real.
Note that when $D' \subset D$, this derivative $\Phi' (0)$ is necessarily smaller than one. In some sense, 
this number is a way to measure ``how smaller $D'$ is compared to $D$, seen from $0$''.

Let us now focus on the set ${\cal G}$ of conformal maps $\Phi$ from the unit disc $\U$ onto a subset $D$ of $\U$ 
 with $0 \in D$, such that $\Phi (0)= 0$ and $\Phi'(0) >0$.
For $a_n = \Phi^{(n)} (0) / n!$, one can write 
$ \Phi (z) = \Phi'(0)z + \sum_{n \ge 2} a_n z^n $
for any $z \in \U$.
Clearly, one can compose two such conformal maps $\Phi_1$ and $\Phi_2$, and 
$\Phi_1 \circ  \Phi_2$ is also a conformal map from $\U$ onto a subset of $\U$ containing the origin.
The set ${\cal C}$ is therefore a semi-group of conformal maps.

Of course, the semi-group ${\cal G}$ is not commutative, and 
it is ``infinite-dimensional''.   
Note that characterizing the set ${\cal G}$ in terms of the coefficients $a_n$ is not a straightforward fact (``can one see for which $a_n$'s, the series is univalent -- i.e. injective --  on the unit disc?''), see for instance the questions related to the Bieberbach conjecture  and its proof by de Branges. This is a typical instance of the difficulty to relate local features of the conformal map (the series expansion) to its global features. 

\subsection {A conformal invariance property}

We are looking for a measure on self-avoiding loops in the plane $\mu$. For any domain open $D$, we define $\mu_D$ as the measure $\mu$ restricted to the set of loops that stay in $D$. We require the following property:

\medbreak
\noindent
{\it $(*)$ For any simply connected set $D$ and any conformal map $\Phi: D \to \Phi (D)$, one has 
$\Phi \circ \mu_D = \mu_{\Phi(D)}$.}
\medbreak

In other words, we want the measures $\mu_D$ to be conformally invariant. As we shall see, this will turn out to be a rather strong constraint.

First of all, the measure $\mu_D$ (and $\mu$ therefore also) must be an infinite measure (if $\mu \not= 0 $). 
This is just because for any $D' \subset D$, the total mass of $\mu_{D'}$ and that of $\mu_D$ must be identical
(because $D$ and $D'$ are conformally invariant). Also, the measure $\mu$ must be scale-invariant and translation-invariant.

Let us now focus on $\mu_\U$. This is now a measure on the set of loops in the unit disc. Consider a map $\Phi 
\in {\cal G}$. $\Phi^{-1}$ is a conformal map from a subset $D$ of $\U$ onto $\U$. One can therefore define the measure $\Phi^{-1} \mu_\U$ as the image of the measure $\mu_D$ under $\Phi^{-1}$. 
Then, the conformal invariance condition is in fact equivalent to the condition that for all $\Phi \in {\cal G}$,
$ \Phi^{-1} \mu_\U = \mu_\U$.
Hence, $(*)$ means that $\mu_\U$ is invariant under the action of the semi-group ${\cal G}$.

\subsection {Uniqueness}

\begin {lemma}
Up to multiplication by a positive real, there exists at most one measure $\mu$ satisfying condition $(*)$.
\end {lemma}

\noindent
{\bf Proof.} We proceed in several steps:
\begin {itemize}
\item
Because of conformal invariance, the whole measure $\mu$ is characterized by $\mu_\U$. Indeed, for any simply connected domain $D$, $\mu_D$ is conformally equivalent to $\mu_\U$, and the knowledge of all $\mu_D$'s characterizes $\mu$.
\item
Because of conformal invariance, the measure $\mu_\U$ is in fact characterized by the measure $\mu_{\U, 0}$ which is defined to be $\mu_\U$ restricted to those loops in $\U$ that disconnect the origin from the boundary of the disc.
Indeed, by conformal invariance, $\mu_{\U, z}$ is conformally equivalent to $\mu_{\U, 0}$ and the knowledge of all 
)$\mu_{\U, z}$'s characterizes $\mu_\U$.
\item
In fact, the knowledge of the quantities  
$$ \mu_{\U, 0} ( \{ w : w \subset \Phi (\U) \} ) $$
where $\Phi \in {\cal G}$ characterizes $\mu_{\U, 0}$ (these events form an intersection-stable family of events that 
generate the $\sigma$-field that one uses).
\end {itemize}
In conclusion, we see that the measure $\mu$ is in fact characterized by the quantities,
$$a(\Phi)=  \mu ( \{w \  : \  w \subset \U,  0 \hbox { is inside } w, w \not\subset \Phi (\U) \})$$
for $\Phi \in {\cal G}$.
Because of conformal invariance again, we see also that 
$$ 
a( \Phi_1 \circ \Phi_2 ) = a (\Phi_1) + a(\Phi_2).
$$
Indeed,
\begin {eqnarray*}
\lefteqn{
\mu ( \{ w \ :  \ w \subset \U , 0 \hbox { is inside } w, w \not\subset \Phi_1 \circ 
\Phi_2 ( \U) \} ) 
}
\\
&=&
\mu ( \{ w \ :  \ w \subset \U , 0 \hbox { is inside } w, w \not\subset \Phi_1  ( \U) \} )
\\
&& 
+ \mu ( \{ w \ :  \ w \subset \Phi_1 (\U) , 0 \hbox { is inside } w, w \not\subset \Phi_1 \circ 
\Phi_2 ( \U) \} ) \\
& = & 
a ( \Phi_1) +
\mu ( \{ w \ :  \ w \subset \U , 0 \hbox { is inside } w, w \not\subset
\Phi_2 ( \U) \} )
\\
&=& a(\Phi_1) + a (\Phi_2) 
\end {eqnarray*}

\begin {figure}[t]
\begin {center}
\begin{picture}(0,0)%
\includegraphics{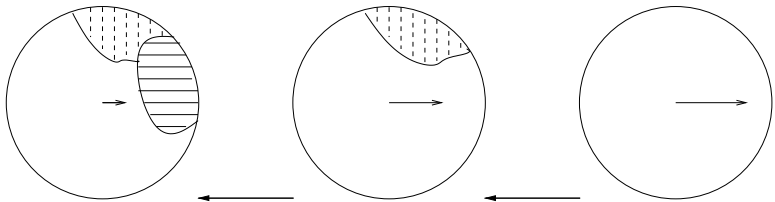}%
\end{picture}%
\setlength{\unitlength}{1973sp}%
\begingroup\makeatletter\ifx\SetFigFont\undefined%
\gdef\SetFigFont#1#2#3#4#5{%
  \reset@font\fontsize{#1}{#2pt}%
  \fontfamily{#3}\fontseries{#4}\fontshape{#5}%
  \selectfont}%
\fi\endgroup%
\begin{picture}(9634,2459)(1184,-2803)
\end{picture}%

\caption {$\Phi_1 \circ \Phi_2$}
\end {center}
\end {figure}

Let us now define the conformal maps 
$\varphi_t  \in {\cal G}$ from the unit disc respectively onto the slit discs
$\U \setminus [1-x_t, 1)$, where $x_t$ is chosen in such a way that 
$\varphi_t' (0)= \exp (-t)$.
Then, clearly, $\varphi_{t+s} = \varphi_t \circ \varphi_s$. This is due to the combination of two facts: Because 
of symmetry, $\varphi_t \circ \varphi_s$ is necessarily one of the $\varphi_u$'s for some $u$, and because of the 
chosen normalization in terms of the derivative at the origin, $u=s+t$.
Hence, it follows easily that for some positive 
constant $c$, $a( \varphi_t ) = ct$.
\begin {figure}[b]
\begin {center}
\begin{picture}(0,0)%
\includegraphics{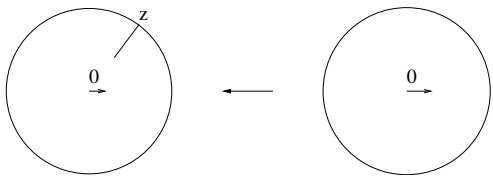}%
\end{picture}%
\setlength{\unitlength}{1381sp}%
\begingroup\makeatletter\ifx\SetFigFont\undefined%
\gdef\SetFigFont#1#2#3#4#5{%
  \reset@font\fontsize{#1}{#2pt}%
  \fontfamily{#3}\fontseries{#4}\fontshape{#5}%
  \selectfont}%
\fi\endgroup%
\begin{picture}(8701,3014)(1508,-2768)
\end{picture}%

\caption {The conformal map $\varphi_{z,t}$}
\end {center}
\end {figure}

Because of rotational invariance, for the same constant $c$, 
$a(\varphi_{z,t}) = ct$ where $\varphi_{z,t} \in {\cal G}$ is this time the conformal map from
$\U$ onto $\U \setminus [(1-x_t) z, z ) $ for a fixed $z \in \partial \U$. 

But, Loewner's theory (we will give more details about this in a second) 
that tells that any conformal map $\Phi \in {\cal G}$ can 
be approximated by iterations of many such maps $\varphi_{z_n,t_n}$.   
It is then not difficult to conclude (as in \cite {LSWsaw} in a slightly different context)
that necessarily $a(\Phi)$ is equal to  $c \log 1/ \Phi'(0)$ for all $\Phi \in {\cal G}$.
Hence, up to this multiplicative constant $c$, there exists at most one candidate 
for a measure $\mu$ satisfying $(*)$. This will conclude the proof of the Lemma. 

\medbreak

It remains to recall a few facts from Loewner's coding of slit domains: Suppose that $\Phi$ is a conformal map in ${\cal G}$ from $\U$ into a set
$\U \setminus \gamma[0,T]$, where $\gamma$ is a smooth simple curve in $\U \setminus \{0\}$ with 
$\gamma(0) \in \partial \U$. 
We can choose the time-parametrization of $\gamma$ in such a way that 
for each $t \le T$, the derivative at zero of the conformal map $\phi_t$  in ${\cal G}$ from $\U$ into  $\U \setminus\gamma [0,t]$ is $\exp (-t)$. In particular, 
$\Phi' (0) = \exp (-T)$.
It turns out to be a little more convenient to work with the inverse maps 
$\psi_t = \phi_t^{-1}$ from $\U$ into $\U \setminus \gamma [0,t]$.

Divide the path $\gamma [0,T]$ into a large number $N$ of portions $\gamma [jT/N, (j+1)T/N]$. 
One can then view 
$\Phi^{-1} = \psi_T$ as the iteration of $N$ conformal maps $\psi^N_1, \ldots, \psi_N^N$, where the map 
$\psi_j^N$ is the conformal map from $\U \setminus \psi_{(j-1)T/N} ( \gamma [(j-1)T/N, jT/N] )$ back onto 
$\U$, such that $\psi_j^N$ maps the origin onto itself and has a positive derivative at the origin.
When $N$ is large, $\psi_{(j-1)T/N} ( \gamma [(j-1)T/N, jT/N] )$ is a small curve in $\U$ starting at 
the boundary point $w(j,N) = w_{(j-1)T/N}= \psi_{(j-1)T/N} (\gamma ((j-1)T/N)))$.
Let us now, for each fixed $N$ and all $j$, replace this small curve by a straight slit 
$[w(j,N) , (1-x_{T/N}) w(j,N)]$, and define the corresponding conformal 
map $\tilde \psi_{j}^N$. And also, define
$$ \tilde \Psi = \tilde \Psi^{T,N} = \tilde \psi_1^N \circ \ldots \circ \tilde \psi_N^N.$$
It turns out that in fact, if one focuses at the behavior of the conformal maps $\psi$ away from the unit circle -- say their restriction to $\U/2$, $\tilde \psi_j^N$ does approximate $\psi_j^N$ up to an error of order $N^{-2}$.
In particular, if one iterates these $N$ terms, one has an error term of the order of $N^{-1}$. 
Hence, when $N$ tends to infinity, the conformal map $\tilde \Psi$ converges (in some suitable sense) towards the 
conformal map $\Psi$. This argument leads to the following two consequences:
\begin {itemize}
\item
$\tilde \Psi' (0) \to \Psi'(0)$ when $N \to \infty$. 
\item
$a ( \tilde \Psi^{-1} ) \to a (\Psi^{-1})$ when $N \to \infty$.
\end {itemize}
But, since $\tilde \Psi$ is the iteration of $N$ conformal maps that remove straight slits, we know that
$
a ( \tilde \Psi^{-1} ) = c \log \Psi' (0).
$.
Hence, $a(\tilde \Phi)$ is indeed equal to $c \log 1/ \Phi'(0)$ for such a slit map $\Phi$.

Finally, one can note that it is possible to approximate any map $\Phi \in {\cal G}$
by such slit maps:
For each $\Phi \in {\cal G}$ and each $\delta > 0$, define the domain 
$D_\delta = \Phi ( (1-\delta) \U )$. It has an analytic boundary, that can be obtained 
via the iteration of slit maps. Furthermore, when $\delta \to 0$, both the 
derivative at the origin and the functional $a$ behave continuously. Hence, 
the formula for $a ( \Phi)$ turns out to be valid for all $\Phi \in {\cal G}$.

\medbreak

We now have to check whether there exists a measure 
$\mu$ for which the formula $a(\Phi) = c \log (1/ \Phi' (0))$ holds.

\subsection {Existence}

{\bf Via Brownian loops:}
Recall that the Brownian loop measure is conformally invariant. This is the scaling limit of the 
measure on random walks loops. It is known since \cite {BL} that a planar Brownian loop has no cut-points.
In particular, the boundary of the unbounded connected component of its complement is a self-avoiding loop in the 
plane. The measure on these loops induced by the Brownian loop-measure is therefore conformally invariant.
Hence, this proves the existence of a measure $\mu$ satisfying $(*)$, since this outer boundary of the Brownian loop
does the job!

Let us briefly give a proper direct definition of the Brownian loop measure  \cite {LWls}:
It will be defined via a product measure. Let $dz$ denote the Lebesgue measure in $\C$, $dT$ the Lebesgue measure
on $\R_+$, and $dP (\gamma)$ the law of a planar Brownian bridge of time-length one starting and ending at the origin
(that is a standard Brownian motion conditioned to be back at the origin at time $1$).  
Note that by scaling and translation, it is possible to transform the standard loop $\gamma$ into a loop $\gamma$ of length T starting and ending at $z$: $\gamma^{z,T} (t) = z+ \sqrt  {T} \gamma (t/T)$.
As before, we then define $\bar \gamma^{z,T}$ from $\gamma^{z,T}$ by identifying two loops if they can be obtained one from the other by time-shift. 

Now, if we work under the product measure $$dz \otimes \frac {dT}{2 \pi T^2} \otimes P(d\gamma),$$ then the measure 
of $\bar \gamma^{z,T}$ is the Brownian loop measure $M$. It will sometimes be convenient to work with multiples
$cM$ of $M$, that we call the Brownian loop measure with intensity $c$.

\begin {figure}
\begin {center}
\includegraphics [width=4in]{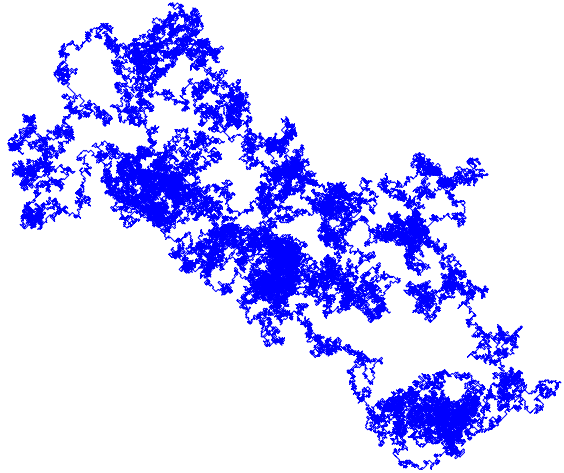}
\caption {A Brownian loop}
\end {center}
\end {figure}

\medbreak
\noindent
{\bf Via percolation clusters:}
Recall the discrete measure on critical percolation clusters $\mu_{pcl}$ on the triangular lattice
and the measure $\mu_{perco}$ on their outer contours. We will be focusing here on the latter.  
When it is defined on the triangular lattice of mesh-size $\delta$,
 we denote it by $\mu_{perco}^\delta$. 
It is a measure on the set of discrete self-avoiding loops.

Camia and Newman \cite {CN} showed how to use  Smirnov's result \cite {Sm}
to deduce results of the following type: The limit $\mu_{perco}^0$ of $\mu_{perco}^\delta$ exists and is 
conformally invariant.
This scaling limit can be in fact described by
 the SLE process with parameter 6 that we will describe in the next section.
It is easy to see that it is not
 supported on the set of self-avoiding loops. Some nearby points
on the discrete scale become double points in the scaling limits. However, the obtained path is 
still a compact loop, and one can define its outer boundary and this one is a self-avoiding loop.

It is in fact not very difficult (using some a priori estimates on five-arm exponents etc.)
to see that the measure $\nu$ on self-avoiding loops that is defined in 
this way is also the scaling limit of discrete ``outer perimeters of clusters''
 obtained by not entering fjords of width equal to the mesh-size of the lattice.  
 
Hence, critical percolation in its scaling limit also defines a measure $\nu$ on self-avoiding loops  
that satisfies $(*)$.

\begin {figure}
\begin {center}
\includegraphics [width=3in]{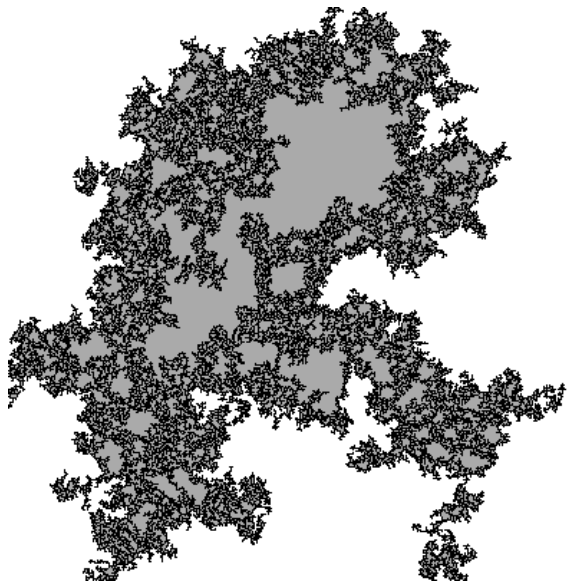}
\caption {A large percolation cluster (pic. by V. Beffara)}
\end {center}
\end {figure}

\medbreak
\noindent
{\bf Consequence.}
It therefore follows from the preceding considerations that:
\begin {itemize}
\item
For some constant $c$, the measure $\nu$ on outer boundaries of percolation clusters is identical to the 
measure of outer boundaries of Brownian loops defined under $cM$. 
\item
If the measure on self-avoiding loops has a scaling limit when the mesh-size of the lattice goes to zero, then 
it is also a multiple of $\mu$.
\end {itemize}
This is probably the simplest and most transparent proof of the equivalence between percolation cluster boundaries, Brownian motion boundaries, and the conjectural scaling limit of self-avoiding walks. This equivalence had been obtained in a different form in \cite {LW2, LSW2, LSWsaw}, for instance to derive \cite {LSW2}
via SLE calculations
(and this is the only existing rigorous proof of this fact at this point) 
the multi-fractal spectrum of the boundary of a Brownian motion with respect to harmonic measure, and in particular 
its Hausdorff dimension $4/3$.
In our setup, it shows that the measure $\mu$ is supported on self-avoiding loops of Hausdorff dimension $4/3$.

\medbreak

Let us also note from \cite {W94} that the whole measure $\mu$ is described by a sample of a single Brownian motion. By looking at the shapes of all connected components of the complement of a sample of a Brownian motion (or a Brownian loop), ordered according to their area (in decreasing order), then the knowledge of all these shapes allows in fact to 
reconstruct the whole measure $\mu$.

\medbreak

Let us now summarize the results of this lecture/section (see \cite {Wsal} for more details):

\begin {theorem}
Up to multiplication by a positive constant, there exists a unique conformally invariant measure on self-avoiding loops in the plane. It can be equivalently described as 
boundaries of Brownian loops, via scaling limits of percolation clusters.
\end {theorem}

A consequence of the study of SLE defined in the next section is in fact that this measure can be defined by SLE means as well, and that it is supported on the set of loops that have Hausdorff dimension $4/3$ (see \cite {LSW4/3}).

We will not go into details here, but this measure and the Brownian loop measure allow to reinterpret in a simple way most nice conformally invariant quantities introduced, such as Schwarzian derivatives etc. and possibly also to introduce new and natural invariants.

Let us however stress an important property of the measure $\mu$:

\begin {theorem}
The measure $\mu$ on self-avoiding loops is invariant under the inversion $z \mapsto 1/z$. More generally, it satisfies the property $(*)$ even for non-simply connected domains $D$.
\end {theorem}

This is not a straightforward result. While the Brownian loop measure $M$ can easily be shown to be 
invariant under the inversion $z \mapsto 1/z$, the ``outside boundary'' of the loop (the boundary of the connected component of its complement that contains infinity) becomes the ``interior boundary'' (the boundary of the connected component of its complement that contains the origin) and the fact that Brownian ``interior boundaries'' and ``outer boundaries'' are the same is rather surprising.
This result has various consequences that are explored in \cite {Wsal}. In particular, it allows to work on any Riemann surface, and to make the theory (loop-soups etc.) work on them.
 
\eject

\section {Schramm-Loewner Evolutions}

\subsection {Definition}

We now give a brief introduction to the SLE process. The reader acquainted to it may want to skip this section.

Suppose that one wishes to describe the law of a random continuous curve $\gamma$ in a simply connected domain $D$, that joins two boundary points $A$ and $B$. Typically, this curve would be the scaling limit of an interface in a discrete grid-based model.
For instance:
\begin {itemize}
\item
For the infinite measures on paths that we have considered (scaling limit of percolation boundaries, Brownian loops, conjectured self-avoiding loops), one way to get such a measure is to start with the discrete measure, and to restrict 
the class of loops to those that contain one chosen part of the boundary of $D$ between $A$ and $B$. In other words, the loop will consist of a path inside $D$ from $A$ to $B$, and it then goes back to $B$ along the boundary of $D$.
Then, renormalize this (finite) measure in order to get a probability measure. So, this path can be viewed as a conditioned version of the loop measure.

\item
In the FK representation of the Potts models, use the following boundary conditions: On one of the two parts of (the discretization of $\partial D$, take ``free'' boundary conditions as before, but identify all vertices on the other part as a single vertex (this is the ``wired'' boundary condition. Then the path that one is interested in in the 
outer contour of the cluster that is attached to the wired part of the boundary. This is indeed a path from $A$ to $B$ in $D$.
\end {itemize}

We are going to assume that this law on continuous curves from $A$ to $B$ in $D$ 
satisfies the following two properties that we shall refer to as 1. and 2.:

\begin {enumerate}
\item 
{\em The law of $\gamma$ is invariant under the (one-dimensional) family of conformal maps from $D$ onto itself that leave $A$ and $B$ unchanged.}

If this property holds, then for any $D', A', B'$ as before, it is possible to define the law of a random curve
that joins $A'$ and $B'$ in $D'$, by taking the conformal image of the law in $D$ (the law of this conformal image is independent of the actual choice of the conformal map because of 1.). 
By construction, the law of this random curve is then conformally invariant in the sense described before (for any triplet $D,A,B$ it defines a law $P_{D,A,B}$ on curves joining $A$ to $B$ in $D$, and for any conformal map $\Phi$,
$\Phi \circ P_{D,A,B} = P_{\Phi(D), \Phi(A), \Phi(B)}$).
Note that it is not really difficult to find a random curve satisfying 1.: Take $D$ to be the upper half-plane $\HH$, $A=0$, $B= \infty$. Then, any random curve $\gamma$ that is scale-invariant in law would do (the conformal maps form $\HH$ onto itself that preserve $0$ and $\infty$ are just the multiplications by a positive constant).

\item 
 {\em The ``Markovian property''}:
Suppose that the curve is parametrized ``from $A$ to $B$'' and  that we know how $\gamma$ begins (i.e. we know $\gamma [0,t]$ say). What is the (conditional) law of the future of $\gamma$?
{\em It is exactly the law $P_{D_t, \gamma_t, B}$, where $D_t$ is the connected component of $D\setminus \gamma[0,t]$ that has $B$ on its boundary}  (when $\gamma$ is a simple curve, then $D_t$ is simply 
$D \setminus \gamma [0,t]$).

\end {enumerate}


It is very natural to combine these two properties. In most of the relevant discrete models that are conjectured to be conformal invariant, one can find natural curves in these models (interfaces etc.) that do satisfy a discrete analog of 2. In particular, critical interfaces in nearest-neighbor interaction models are expected to satisfy both 1. and 2.

Schramm \cite {S1} pointed out that there are not that many possible random curves with these two properties. They form a one-dimensional family and are the so-called SLE curves (SLE stands for Schramm --or Stochastic-- Loewner Evolutions).
These random curves can be described and constructed via iterations of random conformal maps.

Let us now sketch the basic idea behind the construction of SLEs:
Suppose for instance that $D$ is the unit square and that $A$ and $B$ are two opposite corners. For convenience, let us focus on the case where the 
random curve $\gamma$ is simple. Construct first the path $\gamma$ up to some  
time $t$. Define then the unique conformal map $f_t$ from $D \setminus \gamma [0,t]$ onto $D$ such that:
$f_t (\gamma(t))= A$, $f_t (B) =B$ and $f_t'(B) =1$. The image of the path $\gamma [0,t]$ under $f_t$ becomes now a part of the boundary of $D$ (that contains $A$). 

\begin {figure}[t]
\begin {center}
\begin{picture}(0,0)%
\includegraphics{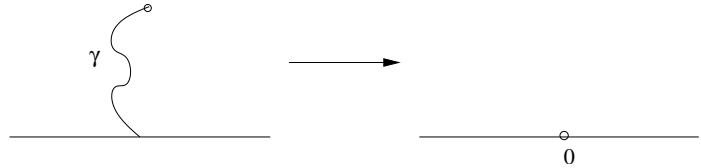}%
\end{picture}%
\setlength{\unitlength}{1539sp}%
\begingroup\makeatletter\ifx\SetFigFont\undefined%
\gdef\SetFigFont#1#2#3#4#5{%
  \reset@font\fontsize{#1}{#2pt}%
  \fontfamily{#3}\fontseries{#4}\fontshape{#5}%
  \selectfont}%
\fi\endgroup%
\begin{picture}(11144,2578)(1179,-2896)
\end{picture}%

\caption {The conformal map $f_t$.}
\end {center}
\end {figure}

Suppose now that we continue the path $\gamma$ after time $t$. Then, the image of this additional stretch under $f_t$ is a simple path in $D$ that joins $f_t (\gamma (t))=A$ to $ f_t (B) = B$. The combinations of the properties 1. and 2.  implies immediately that the conditional law of $\tilde \gamma := f_t ( \gamma [t,\infty))$ (given $\gamma [0,t]$) is $P_{D,A,B}$
(because the conditional law of $\gamma [t,\infty)$ is $P_{D_t, \gamma(t), B}$ and conformal invariance holds). 
In particular, this conditional law is independent of $\gamma [0,t]$. Hence, $\gamma [0,t]$ and $\tilde \gamma [0,t]$ are independent and have the same law.

Now, the conformal map $f_{2t}$ is clearly obtained by composing the two independent identically distributed maps $\tilde f_t$ (corresponding to $\tilde \gamma [0,t]$) and $f_t$. Similarly, for any $N>1$, $f_{Nt}$ is obtained by iterating $N$ independent copies of $f_t$.
It is therefore very natural to encode the curve $\gamma$ via these conformal maps $f_t$ so that the Markovian property translates into independence of the ``increments'' of $f_t$.


In fact, for any $t$, $f_t$ itself can be viewed as the iteration of $N$ independent copies of $f_{t/N}$.
Hence, for any small $\epsilon$, the knowledge of the law of $f_\epsilon$ yields the law of $\gamma$ at any time that is a multiple of $\epsilon$.
This leads to the fact that the knowledge of the law of $f_\epsilon$ for infinitesimal $\epsilon$ in fact 
characterizes the law of the entire curve $\gamma$.

It turns out that there is just ``one possible direction'' in which $f_\epsilon$ can move on infinitesimal level, and this implies that there exists just a one-dimensional family of possible laws for $f_\epsilon$, each one corresponding to a certain speed of motion in this direction. This leads to the one-dimensional family of possible
SLE curves.

More precisely, suppose now that $D$ is the upper half-plane, and that $A=0$ and $B=\infty$. Then, the conformal map
$f_t$ has a Laurent expansion at infinity:
$$f_t (z) = z + a_0 (t) + a_1(t) / z + o (1/z).$$
It is easy to see that $a_1 (t)$ is positive, increasing continuously with $t$, and that it is therefore natural to use a time-parametrization for which $a_1(t)$ increases linearly (in this way, $f_{2t}$ is indeed obtained as the iteration of two independent copies of $f_t$); for historical reasons, the convention is to take 
$a_1(t) = 2t$. Then, the Markovian property implies immediately that $a_0(t)$ is a (real-valued) symmetric Markov process, with independent increments. This yields that $a_0( t) = \beta ( \kappa t)$ for some constant $\kappa \ge 0$ where $\beta$ is a standard real-valued Brownian motion.

Furthermore, one can recover $\gamma$ if one knows only $a_0(t)$: One has to solve (for each fixed $z$) the stochastic differential equation 
\begin {equation}
\label {sleeq}
 df_t (z)  = da_0 (t) + \frac 2 {f_t (z)} dt
 \end {equation}
and see that $\gamma_t = f_t^{-1} (A)$. 
It is not easy but can be proved \cite {RS} that this procedure indeed defines (for any fixed $\kappa$) almost surely a path $\gamma$.

Some of the SLE curves can be shown to have special properties that can be related to the special features of the corresponding discrete models: The independence properties of percolation correspond to the ``locality'' property of SLE$_6$ derived and studied in \cite {LSW1,LSW2}. 
The special properties of self-avoiding walks correspond to the restriction property of SLE$_{8/3}$ related to the properties of self-avoiding walk, e.g. \cite {LSWsaw,LSWrest}.

With this description of the SLE curve via the previous stochastic differential equation, it is not difficult to compute directly the value of critical exponents \cite {LSW1,LSW2,LSW5} that should for instance describe the multi-fractal spectrum of harmonic measure on these curves (and in particular compute their dimension). As a look at the proof in \cite {Be} will probably convince you, there is then still substantial work before one is allowed to state 
rigorously that these numbers are indeed the dimensions of the curves (this last issue is not addressed in the 
theoretical physics papers devoted to this subject where the values of these exponents were first predicted -- see Duplantier's notes \cite {Dqg} and the references therein).

As mentioned in our introduction, we will not go into this part of SLE theory (the interested reader may want to consult \cite {Lbook, Lln, Wln}) and focus on more recent progress. But, we ought to mention here some basic facts about SLE curves (these ``basic'' facts are by no means easy to prove):
The fractal dimension of the SLE$_\kappa$ curve 
(see Rohde-Schramm \cite {RS} for the upper bound and Beffara \cite {Be} for the involved lower bound) is almost surely  $1+ \kappa /8$ when $\kappa < 8$. When $\kappa \ge 8$ the curve becomes space-filling \cite {RS} and has Hausdorff dimension 2. 
This already show that these random curves can have different sometimes surprising fractal structures.
There is also an important phase transition at the value $\kappa=4$:  The SLE curves are simple (without double points) when $\kappa \le 4$, and they have double points as soon as $\kappa > 4$; see \cite {RS}.

\subsection {Computing with SLE}

\begin {figure}
\begin {center}
\begin{picture}(0,0)%
\includegraphics{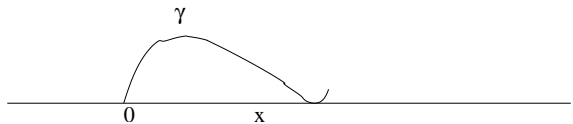}%
\end{picture}%
\setlength{\unitlength}{1539sp}%
\begingroup\makeatletter\ifx\SetFigFont\undefined%
\gdef\SetFigFont#1#2#3#4#5{%
  \reset@font\fontsize{#1}{#2pt}%
  \fontfamily{#3}\fontseries{#4}\fontshape{#5}%
  \selectfont}%
\fi\endgroup%
\begin{picture}(9099,1866)(4714,-4261)
\end{picture}%

\caption {The path disconnects $x$ (sketch)}
\end {center}
\end {figure}

Let us now show why these two phase transitions at $\kappa=4$ and $\kappa=8$ show up. We give a slightly non-rigorous 
description of this derivation in order to be as transparent as possible for those who are not used to 
stochastic analysis, but all the following steps can easily be made rigorous using standard 
stochastic calculus considerations (see \cite {LSW1}). The purpose of the following paragraphs is to 
illustrate how to perform SLE calculations.

The first thing to do is to see what the fact that $f_t (z)$ hits zero at some time $t_0$ (for a fixed $z$) means in terms of the SLE curve $\gamma$. 
Recall that the conformal map $f_t$ maps 
the tip of the curve $\gamma_t$ to $0$, and that it is a conformal map from (the unbounded connected component of) 
$\HH\setminus \gamma[0,t]$ normalized at infinity.
Hence, intuitively, if $f_t(z)$ is close to zero, then it means that ``seen from infinity'' $\gamma_t$ and $z$ are very close in $\HH\setminus \gamma [0,t]$. In fact, what happens is that $f_t (z)$ hits $0$ when the curve $\gamma [0,t]$ disconnects $z$ from infinity in $\HH$. For instance, if $x$ is a positive real, then $f_t(x)$ hits the origin when $\gamma [0,t]$ hits the interval $[x, \infty)$.

Let us now consider two reals $x$ and $y$. $y$ will be positive, and $x <y$ (we will decide later if we choose 
$x$ to be positive or negative). 
Then define 
$$Z_t = \frac {f_t (y) }{ f_t (y) - f_t(x)}.
$$ 
Since $f_0$ is the identity map, $Z_0 = y / (y-x)$. 
It  immediately follows from (\ref {sleeq}) that $d (f_t(y) - f_t (x)) = (2 / f_t (y) - 2 / f_t (x)) dt$
and that 
$$
dZ_t = 
\frac { da_0 (t) } { f_t (y) - f_t (x)} + \frac {2 ( f_t (y) + f_t (x))} {f_t (x) f_t (y) (f_t (y) - f_t (x))    }
.$$
This means that for very small $t$, the law of $Z_t$ can be well-approximated by 
\begin {equation}
\label {taylor}
Z_0 +  \frac {\sqrt {\kappa t} {\cal N}}{y-x} + \frac {2t (y+x)}{ xy (y-x)}
\end {equation}
where ${\cal N}$ is a standard Gaussian random variable with mean $1$ (recall that $a_0(t)$ is a Brownian motion 
with speed $\kappa$).

\begin {figure}
\begin {center}
\begin{picture}(0,0)%
\includegraphics{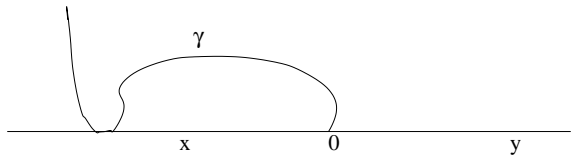}%
\end{picture}%
\setlength{\unitlength}{1539sp}%
\begingroup\makeatletter\ifx\SetFigFont\undefined%
\gdef\SetFigFont#1#2#3#4#5{%
  \reset@font\fontsize{#1}{#2pt}%
  \fontfamily{#3}\fontseries{#4}\fontshape{#5}%
  \selectfont}%
\fi\endgroup%
\begin{picture}(9099,2402)(1414,-4337)
\end{picture}%

\caption {The event $E(x,y)$ when $x<0<y$ (sketch)}
\end {center}
\end {figure}

Now let us consider the event $E(x,y)$ that $f_t (x)$ hits zero strictly before $f_t (y)$ hits zero. This means that the curve $\gamma$ disconnects $x$ from infinity strictly before it disconnects $y$ from infinity.
Because of scale-invariance of the SLE curve, the probability of $E(x,y)$ is a function of $y/x$ only, or equivalently, it is a function of $y / (y-x) = Z_0$. Let us call $G(Z_0)$ this function so that
$$ P ( E (x,y)) = G (y / (y-x)).$$ 

Furthermore, if we see the curve $\gamma$ up to a little time $t$, and look at image of the curve by $f_t$, then we see that 
$E(x,y)$ happens if and only if  $f_t \circ \gamma [t+ \cdot)$ disconnects $f_t (x)$ from infinity strictly before $f_t (y)$.  The Markovian condition then says that this (conditioned) probability is exactly $G(Z_t)$. 
This implies in particular that for small times $t$, $G(Z_0) = E ( G (Z_t))$. 
If we Taylor-expand $G(Z_t)$ at second order using (\ref {taylor}) and 
take the expectation (note that $E({\cal N})=0$ and $E({\cal N}^2) = 1$), we see that 
$$
E(G (Z_t)) \sim G(Z_0)
+  \frac {2(x+y)}{xy(y-x)} G' (Z_0) t + 
\frac {\kappa t } {2 (y-x)^2}  G'' (Z_0)
+ o (t)
.$$
But this quantity is equal to $G(Z_0)$. This implies readily that 
$$ (y-x)^{-2} \left[ \frac {\kappa} 2 G'' (Z_0)  + 2 \left( \frac {1}{Z_0} + \frac {1}{Z_0 - 1} \right) G'' (Z_0) \right] =0 .$$
Since this is true for all possible $Z_0$, we see that $G$ is a solution of the second order 
equation 
$$ \frac {G''(z)}{G' (z)} = \frac {-4}{ \kappa} \left( \frac 1 z + \frac 1 {z-1} \right)
$$ 
i.e. that for some constant $c$, 
\begin {equation}
\label {der}
|G'(z)| = c ( z |z-1| )^{-4/\kappa}.
\end {equation}

Let us now concentrate on the case where $x<0 <y$. 
Here $E(x,y)$ means that the SLE started from the origin hits $(- \infty, x)$ before it does hit $(y, \infty)$.
Of course, this makes sense only if the curve $\gamma$ hits the real axis again at positive times. If this is the case, it is easy to see that 
$P(x,y)$ goes to $0$ when $y \to 0$ (and $x$ is fixed), and to $1$ when $x \to 0$ (and $y$ is fixed).
Hence, $G$ is a function defined on $[0,1]$ satisfying (\ref {der}) with the boundary values $G(0)= 0$ and
$G(1)=1$. 
This means that for $z \in [0,1]$,
$$ G(z)  = c \int_0^z \frac {du}{(u (1-u))^{4/\kappa}}
$$
where $c= \int_0^1 (u (1-u))^{-4 / \kappa} du$. 
This formula makes sense only if $\kappa > 4$ (otherwise the integrals diverge). This corresponds to the fact that 
an SLE with parameter $\kappa$ hits the real axis again only if $\kappa > 4$.  Using the Markovian condition, this turns out to be equivalent to the fact that $\gamma$ is a {\em simple} curve only when $\kappa \le 4$.

\begin {figure}
\begin {center}
\begin{picture}(0,0)%
\includegraphics{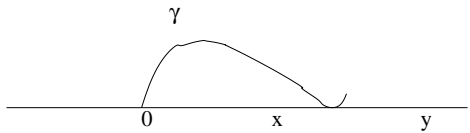}%
\end{picture}%
\setlength{\unitlength}{1539sp}%
\begingroup\makeatletter\ifx\SetFigFont\undefined%
\gdef\SetFigFont#1#2#3#4#5{%
  \reset@font\fontsize{#1}{#2pt}%
  \fontfamily{#3}\fontseries{#4}\fontshape{#5}%
  \selectfont}%
\fi\endgroup%
\begin{picture}(7449,2017)(4414,-4337)
\end{picture}%

\caption {The event $E(x,y)$ when $0<x<y$ (sketch)}
\end {center}
\end {figure}

We now turn our attention to the case where $0 < x < y$. This time, the event $E(x,y)$ means that 
the curve $\gamma$ does hit the interval $(x,y)$ (this clearly is interesting only when $\kappa >4$ as otherwise this probability is $0$). 
The variable $z$ is now in $[1, \infty)$. The boundary condition at $z=1$ is again that $G(1)= 1$. We therefore get that for $z > 1$,
$$
 G(z)  = c \int_z^\infty \frac {du}{(u (u-1))^{4/\kappa} } 
$$ 
where $c = \int_1^\infty (u (u-1))^{-4/ \kappa} du$.
We see that we have this time a problem when $\kappa \ge 8$ because the integral diverges at infinity.
This corresponds to the fact that in this case, the curve is space-filling and hits almost the point $x$ before 
$y$, so that $P(E(x,y))$ is always $1$. Otherwise, the previous formula for $G$ holds.

\medbreak
It is easy to note that
$G(2) = 1/2$ only when $\kappa =6$. Indeed, if $G(2) = 1$, 
then it means that 
$$\int_1^2 \frac {du}{(u (u-1))^{4/\kappa} }=
\int_2^\infty \frac {du}{(u (u-1))^{4/\kappa} }.$$
A simple change of variables (write $v= 1 + (u-1)^{-1}$) shows that the right-hand side
is equal to $\int_1^2 (u (u-1))^{-4/ \kappa} u^{12/ \kappa -2} du$.
Hence, this identity holds if and only if $\kappa =6$.

On the other hand, we know a priori that for  discrete percolation in the upper half-plane in the triangular latticem, the probability of a white crossing from $[0,1]$
to $[2, \infty)$ is $1/2$ (because the complement of this event is the existence of 
a closed crossing from $[1,2]$ to $(-\infty, 0$ which has the same probability by symmetry). Hence, SLE$_6$ is the unique possible scaling limit of critical percolation interfaces. The formula for $G$ in this case is known as Cardy's formula. Smirnov \cite {Sm} has in fact proved that on the triangular lattice, the crossing probabilities indeed converges to $G$ in the scaling limit.
This, together with the construction of the measure $\nu$ in the previous section as outer boundaries of percolation
clusters in the scaling limit, suggest that there is a direct relation between SLE$_6$ and the measure $\nu$ (and this is indeed the case and it can be shown directly using 
SLE$_6$'s locality property and the ideas developped in \cite {LSWrest}). Similarly, SLE$_{8/3}$, which is conjectured to be the scaling limit of self-avoiding walks (see \cite {LSWsaw} for a detailed discussion) is also directly related to the measure $\nu$. 
This is why SLE computations for $\kappa =8/3$ and $\kappa =6$ lead to information 
about the measure $\mu$ on self-avoiding loops (dimension and multifractal spectrum of harmonic measure for these curves etc). 

\eject

\section {Conformal loop-ensembles}

\subsection {Definition}

In the second section/lecture, we were interested in measures supported on the set of (single) loops in the plane, 
motivated by the discrete measures on loops introduced at the beginning. We have identified the unique possible
scaling limit of these measures, under the (sometimes proved) conformal invariance conjecture.
This measure corresponds in the physics literature to models with ``zero central charge'' (this uses 
vocabulary from the closely related representation theory). This means that the probabilistic objects 
that one considers ``do not feel the boundary'' of the domain that they live in.  

We now wish to understand the possible scaling limits of the 
probability measures on configurations for the $O(N)$ models. 
In the present section, we shall define a ``conformal Markov property''
that they should satisfy. Then, we will identify and describe the 
probability measures that do possess this property. This part is based on joint ongoing work with Scott Sheffield \cite {ShW}
(that uses ideas and results from \cite {Wls} on loop-soup clusters and \cite {ScSh} on the Gaussian Free Field).

Let us fix a simply connected domain $D$, say that $D$ is equal to the unit disc $\U$.  
We say that a simple configuration of loops is a family $\Gamma = (\gamma_j, j \in J)$ of self-avoiding 
loops in $D$ such that $J$ is at most countable, the loops are disjoint, and no two loops are nested.
In other words, each point $z \in D$ is either
\begin {itemize}
\item Inside exactly one loop.
\item On exactly one loop.
\item Inside no loop.
\end {itemize}
The term ``simple'' refers to the fact that we require the loops to be non-nested and disjoint.

\begin {figure}[t]
\label {f13}
\begin {center}
\begin{picture}(0,0)%
\includegraphics{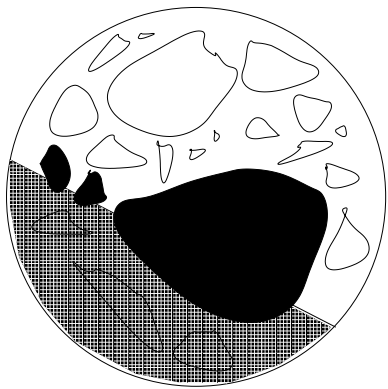}%
\end{picture}%
\setlength{\unitlength}{1539sp}%
\begingroup\makeatletter\ifx\SetFigFont\undefined%
\gdef\SetFigFont#1#2#3#4#5{%
  \reset@font\fontsize{#1}{#2pt}%
  \fontfamily{#3}\fontseries{#4}\fontshape{#5}%
  \selectfont}%
\fi\endgroup%
\begin{picture}(6136,6136)(2633,-5904)
\end{picture}%

\caption {Sketch of the CLE and the domain $\tilde D$.}
\end {center}
\end {figure}

Consider a simply connected subset $D'$ of $D$. For convenience, suppose that $\partial D \cap \partial  D' \not= \emptyset$.
There are two types of loops in $\Gamma$. Those that stay in this smaller domain $ D'$ and those that do exit $D'$. We will call $\tilde J (D')$ and $I(D')$ the corresponding subsets of $J$. In other words:
\begin {eqnarray*}
\tilde J( D', \Gamma) &=& \{ j \in J \ : \ \gamma_j \subset D' \} \\
I (D', \Gamma) &=& \{ j \in J \ : \ \gamma_j \not\subset D' \} 
\end {eqnarray*}
We also define $\tilde D = \tilde D ( \Gamma, D, D')$ the random domain 
that is obtained by removing from $D'$ all  the interiors of 
the loops $\gamma_j, j \in I(D')$.
In other words, for each $\gamma_j$ that does intersect $\partial D'$ (these are the 
loops filled in black in Figure \ref {f13}), 
 we remove all the interior of $\gamma_j$ from $D'$.
Hence, the boundary of $\tilde D$ consists of points 
of $\partial D$ and of points of the loops $\gamma_j, j \in I(D')$, but it contains no point that is inside of a loop.

Note that $\tilde D$ is a (random) domain that can have several connected components 
(for instance if a loop $\gamma_j$ goes in and out of $D'$ several times) and that 
the loops $\gamma_j$ for $j \in \tilde J$ are in the domain $\tilde D$.

\medbreak

Suppose now that $P_\U$ is the law of a random simple configuration of loops
in the unit disc. Assume that $P_\U$ is invariant under any conformal transformation from $\U$ onto itself. Then it is possible to define a law $P_D$ on simple configurations of loops in a simply connected domain $D \not= \C$
by taking the conformal image of $P_\U$.
More generally, for each open domain $D$ in the complex plane, we define 
$P_D$ by taking independent such configurations in each of the connected
components of $D$.

\begin {figure}
\begin {center}
\begin{picture}(0,0)%
\includegraphics{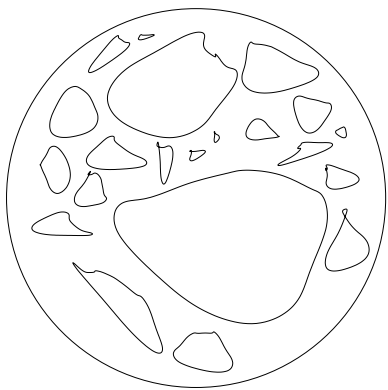}%
\end{picture}%
\setlength{\unitlength}{1539sp}%
\begingroup\makeatletter\ifx\SetFigFont\undefined%
\gdef\SetFigFont#1#2#3#4#5{%
  \reset@font\fontsize{#1}{#2pt}%
  \fontfamily{#3}\fontseries{#4}\fontshape{#5}%
  \selectfont}%
\fi\endgroup%
\begin{picture}(6136,6136)(2633,-5904)
\end{picture}%

\caption {Sketch of a CLE}
\end {center}
\end {figure}

\medbreak
{\bf Definition.}
We say that $\Gamma$ is a simple conformal loop-ensemble if it is a random simple configurations of loops in $\U$
such that:
\begin {itemize}
\item
Its law $P_\U$ is invariant under any conformal map from $\U$ onto itself.
\item
For any $D'$ as above, the conditional law of $(\gamma_j, j \in \tilde J)$ given $(\gamma_j, j \in I(D'))$ is 
$P_{\tilde D}$.
\end {itemize}

\medbreak
It is easy to check that indeed, if the scaling limit of an $O(N)$ model exists and is conformally invariant, then the law of its outermost loops should satisfy this property.
In particular, if one considers the critical Ising model on the triangular lattice, with $1/+$ boundary conditions (i.e. all sites on $\partial D$ are conditioned to have spin $1/+$), then in the scaling limit, the outermost loops surrounding clusters of spin $2/-$
should be a CLE.

For the FK configurations, one would have to concentrate on the critical FK measure with ``free'' boundary conditions, i.e. all edges touching the boundary of the discrete approximation of $D$ are declared to be closed. Then,  assuming conformal invariance,
the scaling limit should satisfy a similar property, but this time, the 
loops are allowed to touch each other (i.e. two loops can have ``touching points'' even if they are not allowed to cross).
We will here restrict ourselves to the ``simple'' case.

\subsection {First properties}

Suppose that $\Gamma$ is such a simple CLE. Suppose also that $\Gamma$ is non-trivial i.e. $P( \Gamma \not= \emptyset 
)>0$. For convenience, we are going to assume that each loop has zero Lebesgue measure.
 
Then:

\begin {itemize}
\item
Almost surely, $\Gamma$ consists of an unbounded countable set of loops.
This follows from the fact that the cardinal of $J$ and the cardinal of $\tilde J$ 
have the same law.

\item
For each fixed point $z \in D$, $z$ is almost surely inside one loop of $\Gamma$.
To see this, let us define the event $A$ that $z$ is surrounded by a loop.
First, since $\Gamma$ is not always empty, conformal invariance implies that one can find 
$D'$ with $z \in D'$, such that with positive probability, there exists a loop surrounding $z$ 
in $I(D')$.
Suppose now that this is not the case, then the conditional probability of $A$ is still $p$ 
because one can start again in $\tilde D$. 
This yields that $P(A)=1$.

\item
The previous statement implies that almost surely, the inside of the loops have full Lebesgue measure, and 
that the Lebesgue measure of the set of points that are outside of all loops is zero.
\end {itemize}

\noindent
Loosely speaking, a CLE defines a fractal-type swiss cheese. 

\medbreak

We are now going to describe how to (try to) define another CLE out of two independent CLEs. 
Suppose that $\Gamma^1$ and $\Gamma^2$ are two independent CLEs (that do not necessarily have the same law). 
Then, $\Gamma^1 \cup \Gamma^2$ consists of a countable family of loops $(\tilde \gamma_j, j \in \tilde J_1 \cup \tilde J_2)$, but
these loops are not disjoint: The loops of $\Gamma^1$ do intersect countably many loops of $\Gamma^2$ and
vice-versa.
We then define the connected components $C_l, l \in L$ of the union of all these loops. Maybe
 there are several disjoint such connected components. In this case, one defines their outer 
boundaries $l_j, j \in L$ (i.e. $l_j$ is the outer boundary of $C_l$). This would be a set of disjoint loops, 
but they might be nested. So we finally keep only the outermost loops of this family in order to have a 
simple configuration of loops. It is then not difficult to see that:

\begin {lemma}
\label {l2}
If these outermost boundaries of unions of independent CLEs form a simple configuration of loops, then it is itself a CLE.
\end {lemma}
 
This lemma shows that in some sense, the set of CLE laws is a semi-group.  
Any two CLEs can be combined into a ``sparser'' CLE (with larger loops). 
This raises for instance the questions of whether a given CLE can itself be viewed as the 
``sum'' of two independent CLEs with the same law. If so, then is any CLE infinitely divisible
(``can it be viewed as an infinite sum of infinitesimal CLEs?'')? 

\subsection {The loop-soup construction}

If one starts to think about what such an infinitesimal CLE might be, it turns out that it should be
related to the infinite measure $\mu$ on self-avoiding loops that we have discussed in an earlier lecture.
This leads to the following construction:

Consider a positive constant $c$, and fix the domain $D= \U$.
Define a Poisson point process of loops with intensity $c \mu_D$ (recall that $\mu_D$ is $\mu$ restricted to those
self-avoiding loops that stay in $D$).
This is a random countable (and overlapping) collection $G$ of self-avoiding loops $(g_k, k \in K)$ in $D$.  
Its law is characterized by the following two facts (if $N(A)$ denotes the 
random number of loops of $G$ that belong to a collection $A$ of loops):
\begin {itemize}
\item
For any two disjoint (and measurable) families of loops $A$ and $B$, $N(A)$ and $N(B)$ are independent.
\item
For any (measurable) family of loops $A$, $N(A)$ is a Poisson random variable with mean $c \mu_D (A)$.
\end {itemize}
It is easy to see (using the additivity properties of Poisson random variables)
 that Poisson point processes exist. Furthermore, if $G_c$ and $G_{c'}$ are two such independent 
Poisson point processes associated to the constants $c$ and $c'$, then $G_c \cup G_{c'}$ is a Poisson point 
process associated to the constant $c_2= c+c'$. In fact, one can construct on the same probability space realizations $G_c$ of the Poisson point processes 
for each $c >0$ in such a way that $G_c$ has intensity $c \mu_D$ and $G_{c_1} \subset G_{c_2}$ for all $c_2 > c_1$.    

One way to think of this is that it rains loops in $D$ with intensity $\mu_D$. $G_c$ is then the family of loops that did fall on $D$ before time $c$. It is called the loop-soup with intensity $c$ in $D$.  It is almost the same as the Brownian loop-soup constructed in \cite {LWls}, except that here, it ``rains'' self-avoiding loops instead of Brownian loops.

\medbreak

When $c$ is large, the set of points in $D$ that is not surrounded by any of the 
loops of $G_c$ is almost surely empty:  This follows rather easily
 from the scale-invariance of the measure $\mu$:
 Take a small $\delta >0$ and define $D_\delta = 
 \{z \in D \ : \  d( z, \partial D) > \delta \}$. 
 For each $n \ge 1$, cover $D_\delta$ by a family ${\cal D}_n$ of $O(4^n)$ discs of radius $\delta 2^{-n}$. 
 For each such disc $d$, one can find at least $n$ concentric disjoint annuli of aspect ratio $2$ in $D \setminus d$.
 Hence, for a positive constant $a$, 
 $$ P ( d \hbox { is not surrounded by a loop in } D ) \le e^{-acn}.$$
 For large $c$, $e^{ac} > 4$.
 Then,
 \begin {eqnarray*}
 \lefteqn {P ( \exists z \in D_\delta \ : \ z \hbox { is not surrounded by a loop in } D ) } \\
 & \le &
 P ( \exists d \in {\cal D}_n \ : \ d  \hbox { is not surrounded by a loop in } D ) \\
 & \le &
 \sum_{d \in {\cal D}_n} 
 P ( d  \hbox { is not surrounded by a loop in } D ) \\
 & \le & O (4^n) \times e^{- acn} 
 .\end {eqnarray*}
 But, if we take $c$ large enough, this last quantity tends to zero as $n \to \infty$. Hence (for large $c$) almost surely, all points in $D_\delta$ (and therefore also in $D$)
 are surrounded by a loop in the loop-soup.
 If we now consider the connected components of $\cup_k g_k$, this implies that for large $c$, there is just one such connected component.

 \medbreak

 On the other hand, when $c$ is very small, things are different:
\begin {lemma}
When $c$ is small (but positive), then 
there almost surely exist a countable family of clusters of loops.
\end {lemma}

In order to prove this, one can couple the loop-soup with the following
 well-known fractal percolation model.
 (In fact, the proof is very close in spirit to that of the corresponding one for 
 multi-scale Poisson percolation in \cite {MR}, Chapter 8.)
 
Consider the unit square $[0,1]^2$. For each $n$, we will divide it into $4^n$ squares 
of side-length $2^{-n}$. To each such square $C$, associate a Bernoulli random variable $X(C)$ equal to 
one with probability $p$. We assume that all $X(C)$'s are independent.
Then, define 
$$
 M= [0,1]^2 \setminus \bigcup_{C \ : \ X(C)=0} C .
$$

\begin {figure}
\begin {center}
\includegraphics [width=2in]{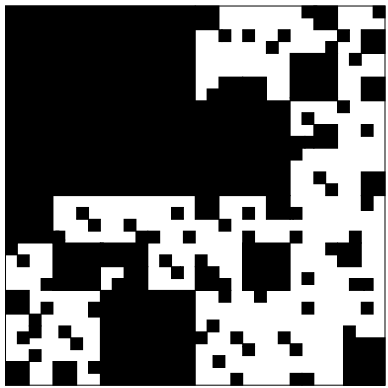}
\caption {The set $M$ (sketch).}
\end {center}
\end {figure}

This is the fractal percolation model introduced by Mandelbrot in \cite {Ma}.
It is very easy to see that the area of $M$ is almost surely zero as soon as $p>0$. 
Chayes, Chayes and Durrett \cite {CCD} 
have shown that this model exhibits a first-order phase transition:
There exists a $p_c$, such that for all $p \ge p_c$, $M$ connects two opposite sides of the 
unit square with positive probability, whereas for all $p < p_c$, the set is a.s. disconnected
(note that in fact, if $p \le 1/4$, then $M$ is almost surely empty -- one can just compare the process with a 
subcritical branching process). Note that at the critical probability, $M$ still has non-trivial connected 
components.

Now, let us consider a loop-soup with intensity $c$ in the unit square. This is just the conformal image of the 
loop-soup in $\U$.
For each loop $g$, let $d(g)$ denote its diameter, and define $n(g)$ in such a way 
that $d(g) \in [2^{-n(g)-1}, 2^{-n(g)})$.
Note that $g$ can intersect at most 4 different dyadic squares with side-length $2^{-n}$. 
We each dyadic square $C$ of side-length $2^{-n}$, we now define 
$$
 \tilde X (C) = 1 - 1_{\exists g\ : \ n(g)= n \hbox { and } C \cap g \not= \emptyset }.
$$ 
The variables $\tilde X (C)$ are of course 
not independent: When $C$ and $C'$ are adjacent and have the same side-length
$2^{-n}$, then 
it can happen $\tilde X(C)= \tilde X (C') = 0$ because a loop $g$ with $n(g) = n$ intersects them both.
But it is easy to see that they are positively correlated. In particular, this implies that it is possible
to couple 
the family $(\tilde X(C))$ with a family of independent Bernoulli random variables $X(C)$ such that for each $C$,
$\tilde X(C) \ge X(C)$ almost surely, and $P( \tilde X (C) = 1) = P ( X(C) = 1)$.
Because of the scale-invariance of the Poisson point process, we see that for some constant $b$, 
$$
P ( \tilde X (C) = 1 ) \ge \exp (-bc ) 
$$ 
for all $C$ (the inequality is just due to the squares adjacent to the boundary of the unit square that play aslightly different role).
 
\begin {figure}
\begin {center}
\includegraphics [width=2in]{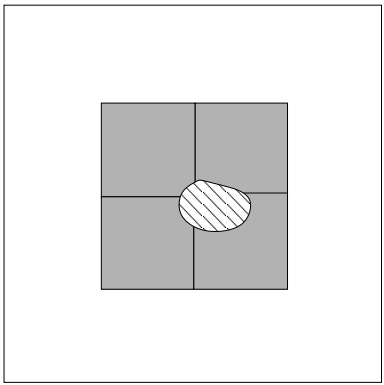}
\caption {A loop and its four corresponding squares (sketch)}
\end {center}
\end {figure}
Hence, we see that the loop-soup percolation process is dominated by a Mandelbrot-percolation 
model (remove each square independently with probability $1-p = 1-\exp (-bc)$). 
In other words, if instead of removing the loops, we remove the squares of the corresponding side-length that they 
intersect, we wtill remove less squares than in the corresponding fractal percolation
model. When $p \ge p_c$, we know that the remaining domains contains paths and loops with positive probability.
Hence, so does  the complement of the filled loops, as soon as 
$\exp (-bc) \ge p_c$ (and this is certainly the case for very small $c$).

One can then upgrade this argument to conclude that
this complement in fact {\em almost surely} contains paths, and loops.

\medbreak

Note that this argument in fact shows that for any $\eps >0$, one can find $c$ small enough such that the 
probability that (if one considers the loop-soup of intensity $c$ in the unit square), there exists a separating path 
from the top to the bottom of the square, that stays in the middle-strip $[1/3,2/3] \times [0,1]$
is greater than $1 - \eps$.
This can then be used to prove that for small $c$, 
all clusters of loops are a.s. at positive distance 
from each other and from the boundary.
\medbreak

If we now consider a domain $D' \subset D$, and a loop-soup in $D$ for such a small $c$. The family of loop-soup cluster outer boundaries are denoted by $\Gamma =(\gamma_j, j \in J)$. Then, for each open subset $D_0$ of $D'$, we see that 
conditionally on $D_0 \subset D_\Gamma'$, the law of the loops $g_k$ that stay in $D_0$ is just that of a loop-soup in $D_0$ (this is because of the Poissonian nature of the loop-soup).
Since this is true for all $D_0$ and using the conformal invariance of the measure $\mu_D$, one can deduce that in fact, the family $\Gamma$ is a CLE.

\begin {theorem}
\label {t5}
There exists $c_0 > 0$ such that:
\begin {itemize}
\item For all $c > c_0$, there exists just one loop-soup cluster.
\item For each $c \le c_0$, the joint law of the outer boundaries of loop-soup clusters defines a CLE.
\end {itemize}
\end {theorem}

One just has to prove that at the critical value $c_0$, one still has several loop-soup clusters. This can be proved in an analogous way as the corresponding result for Mandelbrot percolation in \cite {MR}.

\medbreak

Note that if $c_1$ and $c_2$ are chosen in such a way that $c_1 + c_2 \le c_0$. Then, the CLE defined for the loop-soup 
with intensity $c_1+ c_2$ can be defined starting from two independent CLEs corresponding to loop-soups with intensity $c_1$ and $c_2$ 
just as defined in Lemma \ref {l2}. Hence, the conditions of this lemma are sometimes true...

\medbreak
This loop-soup construction shows that CLEs do exist. It is of course natural to ask whether there exist other CLEs. Also, one might wish to 
study the geometry of the loops in a CLE.

Note that the scaling limit of $O(N)$ models are not simple families of loops (one allows nested ones). However, 
if the scaling limit of the outermost loops in a domain exist and are conformally invariant, and are therefore a CLE, then one can construct the scaling limits of the inner ones by iteratively define ``independent'' CLEs inside each of the outermost loops and so on. 
This defines a ``nested CLE''.

\subsection {The Gaussian Free Field}

The following is a very heuristic description of the rich 
ongoing work by Oded Schramm and Scott Sheffield, and by Scott Sheffield alone on the relation between SLE and the Gaussian Free Field. More details will be found in their upcoming papers \cite {ScSh, Sh}.

\subsubsection {Definition}

Consider a simply connected domain $D$ in the plane such that $D \not= \C$. 
The Gaussian Free Field $\Gamma$ in $D$ (with Dirichlet boundary condition) is a Gaussian 
process (i.e. a Gaussian vector in an infinite-dimensional space) that can be viewed as a 
random generalized function $h$ in $D$. It is not really a function because $h$ is not 
defined at points, but for each smooth compactly supported function $f$ in $D$,
$\Gamma (f)$ is a centered Gaussian random variable, with 
covariance
$$
E ( \Gamma (f)^2 ) = \int_{D \times D} d^2 x  d^2 y f(x) f(y) G_D(x,y)
$$
where $G_D(x,y)$ denotes the Green's function in $D$. 
Since, $G_D(x,x)= \infty$, one can not make sense of $\Gamma (\delta_x)$ (which would then be
$h(x)$). However, it is often useful to think of $\Gamma (f)$ as a generalized Gaussian function $h$ 
with covariance given by $G_D$, so that
$$\Gamma (f) = \int_D h(x) f(x) d^2x .$$
The Gaussian Free field (GFF in short) inherits conformal invariance properties from
the conformal invariance of the Green's function.

One equivalent way to define things is to consider the Hilbert space ${\cal H}_D$ 
of generalized functions in $D$ with zero boundary conditions endowed with the norm $$(h,h)_D = \int_D |\nabla h(x)|^2 d^2 x$$
 (i.e. the completion of some space of smooth functions 
with respect to this norm). Then, $\Gamma$ is simply a Gaussian process in this 
space such that the covariance of $h \in {\cal H}_D$ is $(h,h)_D$ (see \cite {Shgff}).
Then, clearly, the energy is conformally invariant i.e. 
$(h \circ \Phi^{-1} , h \circ \Phi^{-1} )_{\Phi(D)} = (h, h )_D$. 

\medbreak

One important additional property of the GFF is its ``Markovian property''.
A basic version of this property follows rather directly from the 
properties of Gaussian processes (and the fact that in this Gaussian spaces setup, independence and orthogonality are the same).
Suppose for instance that $\gamma$ denotes a smooth deterministic loop (or a simple curve joining two boundary points) in $D$. Let $D_1$ and $D_2$ denote the two 
connected components of $D \setminus \gamma$.
One can decompose the space ${\cal H}_D$ into three orthogonal subspaces of generalized functions:
\begin {itemize}
\item
The spaces ${\cal H}_{D_1}$ and ${\cal H}_{D_2}$ of (generalized) functions in $D_1$ and $D_2$ with Dirichlet boundary conditions 
on their boundaries.
\item
The space ${\cal H}_\gamma$ of (generalized) functions in $D$, that are harmonic in $D_1$, in $D_2$ and have Dirichlet boundary conditions on $\partial D$.
\end {itemize}
In other words, one decomposes $h$ into its value on $\gamma$ (extended in a harmonic way on the complement of $\gamma$)
and then the two contributions in $D_1$ and $D_2$.
For a Gaussian Free Field, these three generalized distributions are 
therefore independent.
Hence,  if one conditions on the values $h_\gamma$ of $h$ on (an infinitesimal neighborhood of) $\gamma$, 
then the field in $D_1$ and $D_2$ are conditionally independent and obtained by taking the sums of a GFF in $D_1$ (resp. $D_2$) with the harmonic extension of 
$h_\gamma$.

Hence, the GFF possesses both a conformal invariance property and a Markovian type property. 
If one thinks of it in terms of a  generalized surface, if one could define some ``geometric'' curves/loops, like level lines, these would be candidates for SLEs/CLEs. 

\subsubsection {``Cliffs'' as level lines}

Suppose that $D=\HH$. Then, consider the field $F$ obtained by taking the sum of a GFF in $\HH$ with the deterministic harmonic function $F_0$ on $\HH$ that takes the boundary value $\lambda$ on $\R_+$ and $0$ on $\R_-$ i.e.
$F_0 (z) = \frac {\lambda}\pi \arg (z)$.
When $\eta$ is a simple curve from $0$ to $\infty$ in $\HH$, then we call $H_-$ and $H_+$ the two connected components
of $\HH \setminus \eta$ (in such a way that $\R_+ \subset \partial H_+$). Then:

\begin {theorem}[\cite {ScSh}]
Almost surely, for one well-chosen value of $\lambda$,
 there exists a unique simple curve $\eta$ (which is a measurable function of $F$)
from the origin to $\infty$ in $\HH$  
such that 
the respective conditional laws (given $\gamma$) of $F$ and $F-\lambda$ 
in the two connected components $H_-$ and $H_+$
are that of independent Gaussian free fields $F_-$ and $F_+$ with Dirichlet boundary conditions. 
 Furthermore, the law of $\eta$ is that of an SLE with parameter $4$.
\end {theorem}

In other words, if one chooses the height-gap $\lambda$ well, then one can follow a ``cliff'' with height-gap $\lambda$
from the origin to infinity.    
The proof of this result is far from obvious. An important step is to make sense of the 
Markovian property for ``random curves that are defined by the field themselves'' that Schramm and Sheffield call 
local sets of the Gaussian Free Field (instead of deterministic curves $\gamma$ as in the last subsection).

It is however possible to make a little calculation that explains why SLE with parameter $\kappa =4$ shows up.
Indeed, recall the SLE equation (\ref {sleeq}). 
If we formally
expand $d (\log f_t (z))$ taking into account the second order terms due to the quadratic variation of $a_0$,
we get that 
$$ d ( \log f_t (z)) = \frac { df_t (z) } { f_t (z)} - \frac {\kappa}{2 f_t^2(z)} dt 
= \frac {da_0 (t)}{f_t (z)} + \frac {2 - \kappa/2}{f_t^2 (z)} dt 
$$
so that $\log f_t (z)$ is a martingale only when $\kappa = 4$. Such a property has to be true along the cliff
$\eta$ because the expected height of $z$ given $\eta [0,t]$ would then be proportional to $\arg (f_t (z))$.

This theorem shows that it is indeed possible to construct geometrically defined curves in a GFF. In fact, it turns out that one can not only define such ``flat cliffs'' but also cliffs with varying height (depending on the winding of the curve) that we shall not describe here. Hence, a realization of the GFF contains (simultaneously) a huge quantity of such geometrically defined curves.

It is therefore natural to try to find CLEs in a GFF. This indeed turns out to be possible. The simplest ones are the ``loop versions'' of the curves defined in the last theorem. Intuitively, if one draws all the loops
in a GFF corresponding to cliffs between height $0$ (on the outside of the loop) and $\pm \lambda$ (on the inside of the loop), the Markovian property of the GFF should imply that the obtained family of loops is a CLE. 

In \cite {ShW}, we plan to write up the proof of the fact that the simple CLEs that one can define in a GFF are the same (in law) as those defined via loop-soups in the previous section. This will show that the critical value of $c$ in Theorem \ref {t5}
as $c_0=1$, and describe the outer boundaries of CLEs as SLE-type loops with parameter 
$\kappa = \kappa (c) \in (8/3, 4]$ where $c= (\kappa - 8/3)(6- \kappa) / ( 2 \kappa)$.
Furthermore, the CLEs defined in the previous section are the only simple CLEs, and therefore the only possible conformally invariant scaling limits of the $O(N)$ models.

\eject

----------------------------

w.w.

Postal address:

D\'epartement de Math\'ematiques

Universit\'e Paris-Sud 

91405 Orsay cedex, France 

\medbreak

wendelin.werner@math.u-psud.fr

\begin {thebibliography}{99}

\bibitem {Be}
{V. Beffara (2002),
the dimensions of SLE curves, preprint.}

\bibitem{BPZ}
{A.A. Belavin, A.M. Polyakov, A.B. Zamolodchikov (1984),
Infinite conformal symmetry in two-dimensional quantum field theory.
Nuclear Phys. B {\bf 241}, 333--380.}

\bibitem {BL}
{K. Burdzy, G.F. Lawler (1990),
Non-intersection exponents for random walk and Brownian motion II. Estimates and 
application to a random fractal, Ann. Prob. {\bf 18}, 981-1009.}

\bibitem {CN}
{F. Camia, C. Newman (2005),
 The Full Scaling Limit of Two-Dimensional Critical Percolation, preprint.}
 
\bibitem {Ca0}
{J.L. Cardy (1984),
Conformal invariance and surface critical behavior,
Nucl. Phys. {\bf B 240}, 514-532.}

\bibitem {Cln}
{J.L. Cardy (2004),
SLE for theoretical physicists, Ann. Physics, to appear.}

\bibitem {CCD}
{J.T. Chayes, L. Chayes, R. Durrett (1988),
Connectivity properties of Mandelbrot's percolation process,
{Probab. Theory Related Fields} {\bf  77},   307-324.}

\bibitem {Dqg}
{B. Duplantier (2004),
 Conformal fractal geometry and boundary quantum gravity, 
  in Fractal Geometry and applications, a jubilee of Beno\^\i t Mandelbrot, 
Proc. Symp. Pure Math. {\bf 72}, vol. II, AMS.}

\bibitem {GT}
{C. Garban, J.A. Trujillo-Ferreras (2005),
The expected area of the Brownian loop is $\pi/5$,
preprint.}

\bibitem {G}
{G.R. Grimmett, Percolation, Springer, 1989.}

\bibitem {G2}
{G.R. Grimmett (1997),
Percolation and disordered systems,
Ecole d'\'et\'e de Probabilit\'es de St-Flour XXVI, L.N. Math. {\bf 1665},
153-300}
 
\bibitem {KNln}
{W. Kager, B. Nienhuis (2004),
A guide to stochastic L\"owner evolution and its applications,
J. Stat. Phys. {\bf 115}, 1149-1229}

\bibitem {Kbook}
{H. Kesten (1984), 
Percolation theory for mathematicians, Birkha\"user.}

\bibitem {Lln}
{G.F. Lawler (2004), 
An introduction to the stochastic Loewner evolution,
in Random Walks and Geometry, 263-293, de Gruyter, Berlin.}

\bibitem {Lbook}
{G.F. Lawler (2005),
Conformally invariant processes in the plane,
Math. Surveys and Monographs {\bf 144}, AMS.}

\bibitem {LSW1}
{G.F. Lawler, O. Schramm, W. Werner (2001),
Values of Brownian intersection exponents II: Half-plane exponents,
Acta Mathematica {\bf 187}, 236-273.}

\bibitem {LSW2}
{G.F. Lawler, O. Schramm, W. Werner (2001),
Values of Brownian intersection exponents II: Plane exponents,
Acta Mathematica {\bf 187}, 275-308.}

\bibitem {LSW4/3}
{G.F. Lawler, O. Schramm, W. Werner (2001),
The dimension of the Brownian frontier is $4/3$,
Math. Res. Lett. {\bf 8}, 401-411.}


\bibitem {LSW5}
{G.F. Lawler, O. Schramm, W. Werner (2002),
One-arm exponent for critical 2D percolation,
Electronic J. Probab. {\bf 7}, paper no.2.}

\bibitem {LSWlesl}
{G.F. Lawler, O. Schramm, W. Werner (2004),
Conformal invariance of planar loop-erased random
walks and uniform spanning trees, Ann. Prob. {\bf 32}, 939-996.}

\bibitem {LSWsaw}
{G.F. Lawler, O. Schramm, W. Werner (2004),
On the scaling limit of planar self-avoiding walks, 
 in Fractal Geometry and applications, a jubilee of Beno\^\i t Mandelbrot, 
Proc. Symp. Pure Math. {\bf 72}, vol. II, AMS 339-364.
}

\bibitem {LSWrest}
{G.F. Lawler, O. Schramm, W. Werner (2003),
Conformal restriction properties. The chordal case,
J. Amer. Math. Soc., {\bf 16}, 917-955.}

\bibitem {LT}
{G.F. Lawler, J.A. Trujillo-Ferreras (2005),
Random walk loop-soup, Trans. A.M.S., to appear.}

\bibitem {LW1}
{G.F. Lawler, W. Werner (1999),
Intersection exponents for planar Brownian motion,
Ann. Probab. {\bf 27}, 1601-1642.}

\bibitem {LW2}
{G.F. Lawler, W. Werner (2000),
Universality for conformally invariant intersection
exponents, J. Europ. Math. Soc. {\bf 2}, 291-328.}

\bibitem {LWls}
{G.F. Lawler, W. Werner (2004),
The Brownian loop-soup, 
Probab. Th. Rel. Fields {\bf 128}, 565-588.}


\bibitem {Ma}
{B.B. Mandelbrot,
{\em The Fractal Geometry of Nature},
Freeman, 1982.}

\bibitem {MR}
{R. Meester, R. Roy,
Continuum Percolation,
CUP, 1996.}
 
\bibitem {N}
{B. Nienhuis (1982),
Exact critical exponents for the $O(n)$ models in two dimensions,
Phys. Rev. Lett. {\bf 49}, 1062-1065.}

\bibitem {RS}
{S. Rohde, O. Schramm (2005), 
Basic properties of SLE, Ann. Math. {\bf 161}, 879-920.}

\bibitem {S1}{
O. Schramm (2000), Scaling limits of loop-erased random walks and
uniform spanning trees, Israel J. Math. {\bf 118}, 221-288.}

\bibitem {Shgff}
{S. Sheffield (2003),
Gaussian Free Fields for mathematicians, preprint.}


\bibitem {ScSh}
{O. Schramm, S. Sheffield (2005), in preparation}

\bibitem {Sh}
{S. Sheffield (2005), in preparation}

\bibitem {ShW}
{S. Sheffield, W. Werner (2005), in preparation}

\bibitem {Sm}
{S. Smirnov (2001),
Critical percolation in the plane: conformal invariance,
 Cardy's formula, scaling limits,
 C. R. Acad. Sci. Paris S�. I Math. {\bf 333},  239-244.}

\bibitem {W94}
{W. Werner (1994),
Sur la forme des composantes connexes du com\-pl\'ementaire de la courbe 
brownienne plane, 
Probab. Th. rel. Fields {\bf 98}, 307-337.}

\bibitem {Wln}
{W. Werner (2004),
Random planar curves and Schramm-Loewner Evolutions,
in 2002 St-Flour summer school,  L.N. Math. {\bf 1840}, 107-195.}

\bibitem {Wls}
{W. Werner (2003),
SLEs as boundaries of clusters of Brownian loops, 
C.R. Acad. Sci. Paris, Ser. I Math. {\bf 337}, 481-486.}

\bibitem {Wln2}
{W. Werner (2005),
Conformal restriction and related questions, 
Probability Surveys {\bf 2}, 145-190.}

\bibitem {Wsal}
{W. Werner (2005),
The conformal invariant measure on self-avoiding loops, 
preprint.
}


\end{thebibliography}
\end{document}